\documentclass[12pt]{article}
\usepackage[T2A]{fontenc}
\usepackage[utf8]{inputenc}
\usepackage{amssymb,amsmath,amsfonts,amsthm,amscd,latexsym,verbatim,graphics,epsfig,indentfirst,xcolor,array}
\usepackage{geometry}
\def\Aut{\mathrm{Aut}}

\def\Com{\mathrm{Com}}

\begin{document}

\sloppy


\begin{center}
{\Large Generalized sharped cubic form and split spin factor algebra}

\smallskip

Vsevolod Gubarev, Farukh Mashurov, Alexander Panasenko
\end{center}

\begin{abstract}
There is a well-known construction of a Jordan algebra via a sharped cubic form.
We introduce a generalized sharped cubic form and prove that the split spin factor algebra 
is induced by this construction and satisfies the identity
$((a,b,c),d,b) + ((c,b,d),a,b) + ((d,b,a),c,b) = 0$.
The split spin factor algebras have recently appeared in the classification of 2-generated 
axial algebras of Monster type fulfilled by T. Yabe; their properties were studied by J. McInroy and S. Shpectorov. 

{\it Keywords}:
sharped cubic form, split spin factor algebra, Lie triple system.
\end{abstract}

\section{Introduction}

The structure theory of nonassociative algebras and rings has been generally based on consideration 
of a concrete variety~$\mathcal{M}$ of algebras and then on the description or a~search of (finite-dimensional) 
simple algebras from~$\mathcal{M}$.
We may say that the structure theory of Jordan, alternative, Malcev etc. algebras was developed more or less in such way.

Recently, a new approach to get plenty of simple nonassociative algebras appears. 
The notion of axial algebra was proposed by J.I. Hall, F. Rehren and S. Shpectorov in 2015~\cite{HRS}. 
It defines a class of (non-associative) commutative algebras generated by specific idempotents, 
and the product in an axial algebra satisfies some restrictions depending on its type. 
Roughly speaking, we may say that these restrictions generalize the ones originated from 
Pierce decomposition fulfilled on associative, alternative, or Jordan algebras.
Axial algebras of Jordan type are close to Jordan algebras, while axial algebras of Monster type 
generalize the Griess algebra.
The latter has the Monster group exactly as its automorphism group.

In the direction of axial algebras, different classifications of algebras with 
a~small number of generators were stated. 
One of them, obtained by T. Yabe in 2020 (published in 2023~\cite{Yabe}), 
provides a list of all 2-generated axial algebras of Monster type $(\xi,\eta)$ 
admitting a flip between generating axes. One of the algebras from the Yabe's list 
was denoted as $S(\alpha,E)$ and the properties of this algebra were studied 
by J.~McInroy and S.~Shpectorov in 2022~\cite{McInroy}.
The authors called them as split spin factor algebras by analogy with spin factor algebra, 
the simple Jordan algebras of special form.

The main goal of the current work to study the identities fulfilled 
on the split spin factors $S(\alpha,t,E)$, where $E$ is any vector space of dimension 
at least two endowed with a symmetric nondegenerate form $\langle\cdot,\cdot\rangle$ 
and $\alpha,t$ are parameters from the ground field~$F$.
We show that there are no identities of degree 3 and 4 on $S(\alpha,t,E)$, $\alpha,t\notin\{0,1\}$, 
which do not follow from commutativity. 
Further, we prove that all identities on~$S(\alpha,E)$, $\alpha\notin\{-1,0,1/2,1,2\}$, 
of degree~5 follows from commutativity and
$$
((a,b,c),d,b) + ((c,b,d),a,b) + ((d,b,a),c,b) = 0,
$$
where $(a,b,c) = (ab)c - a(bc)$. 
We name it as the three associators identity.

In~1965, J.M. Osborn gave the list of all irreducible relative to commutativity identities 
of degree 5~\cite{Osborn}. The fourth of these five identities~\cite[eq.\,(15)]{Osborn} with 
$\delta_2 = -\delta_1\neq0$ is the three associators identity 
with one of the three variables $a,c,d$ equal to $b$, e.\,g., $d = b$.

There is a construction of a Jordan algebra by any sharped cubic form~$(N,\#,c)$~\cite{McCrimmon2}. 
For the proof, you need to verify dozens of relations in terms of $N$, its derived maps $T,S$, $\#$, 
the triple product $\{\cdot,\cdot,\}$ and the $U$-operator, see~\cite[Appendix C]{McCrimmon}.
This construction, in particular, allows to build simple Jordan algebras of Albert type. 
In the work, we consider generalized sharped cubic form and prove the analogues of the known 
relations, which hold for sharped cubic forms. 
With the help of these relations (Lemmas 1--7), we more or less reduce checking 
the three associators identity on the algebra $S(\alpha,t,E)$ to the properties 
of the Psi-map~$\Psi(r,s,q)$, which is defined via an associator, see~\S5. 
Further, in Lemma~9, we show that $\Psi(r,s,q)$ may be expressed in terms of the projections 
of $r,s,q$ on~$E$ and its bilinear form $\langle \cdot,\cdot\rangle$. 
Hence, $\Psi(r,s,q)$ defines a Lie triple product on $S(\alpha,t,E)$ and 
it is connected with the simple pre-Lie algebra from~\cite{Svinolupov}. 
Finally, we prove in Theorem~3 that the three associators identity holds on~$S(\alpha,t,E)$.

Let us provide a short outline of the work.
In~\S2, we give the definition of the split spin factor algebra~$S(\alpha,E)$ 
and its natural generalization~$S(\alpha,t,E)$.
In~\S3, we recall the results about sharped cubic form and induced Jordan algebra.
In~\S4, we introduce a~generalized sharped cubic form and prove the main relations devoted to it. 
In~\S5, we define the Psi-map and derive the equalities on it. 
The goal of~\S6 is to prove that the three associators identity holds on~$S(\alpha,t,E)$. 
In~\S7, with the help of computer algebra, we prove that all identities of degree not greater than 5 
fulfilled on~$S(\alpha,E)$ follow from commutativity and the three associators identity. 
In the general case of~$S(\alpha,t,E)$, there are identities of degree~5, 
which do not follow from commutativity and the three associators identity.
In~\S8, we formulate open problems concerning generalized sharped cubic form and induced algebras. 
In particular, the following question remains to be open: what identity holds on all algebras 
induced by a generalized sharped cubic form? 

We assume that the ground field~$F$ is of characteristic not~2 and not~3.

\section{Split spin factor algebra}

Let~$F$ be quadratically closed, i.\,e. roots of any quadratic equation over~$F$ lie in it.

Below, we put the multiplication table for $S(\alpha,E) = Fz_1+Fz_2+E$, where $\dim E = 2$
and there exists a~nondegenerate bilinear form $\langle\cdot,\cdot\rangle$ on $E$:
\begin{equation} \label{Split-Spin-Prod-Def}
\begin{gathered}
z_1^2 = z_1, \quad
z_2^2 = z_2, \quad
z_1z_2 = 0, \quad
ez_1 = \alpha e, \quad
ez_2 = (1-\alpha)e, \\
ef = -\langle e,f\rangle(\alpha(\alpha-2)z_1+(\alpha^2-1)z_2),\ e,f\in E.
\end{gathered}
\end{equation}

For $\alpha\neq-1,2$, the algebra~$S(\alpha,E)$ admits the nondegenerate invariant bilinear form given by
\begin{equation} \label{Split-Spin-Inv-Form}
\begin{gathered}
(z_1,z_1) = \alpha+1, \quad
(z_2,z_2) = 2-\alpha, \quad
(z_1,z_2) = 0, \\
(e,f) = (\alpha+1)(2-\alpha)\langle e,f\rangle, \quad
(e,z_i)= 0, \ e,f\in E.
\end{gathered}
\end{equation}
Invariancy means that $(ab,c) = (a,bc)$ for all $a,b,c\in S(\alpha,E)$.

Note that in~\eqref{Split-Spin-Prod-Def} and~\eqref{Split-Spin-Inv-Form} 
we may consider any vector space~$E$, not necessarily of dimension~2.
More generally, we may study the algebra~$S(\alpha,t,E)$ depending on two parameters~$\alpha,t$, 
where $E$ is a vector space of any dimension endowed with a nondegenerate bilinear form.
Then the product of elements $e,f\in E$ is defined by the formula
$$
ef = \langle e,f\rangle(z_1+tz_2).
$$
The split spin factor algebra~$S(\alpha,E)$ is an algebra $S(\alpha,t,E)$ 
with $t = (\alpha^2-1)/\alpha(\alpha-2)$.
Surely, we require that $\alpha\neq0,2$.

{\bf Proposition 1}. 
Let $E$ has a finite dimension $n\ge 1$. Then $S(\alpha,t,E)$ is a simple algebra 
if and only if $\alpha\notin \{0,1\}$, $t\neq 0$.

{\sc Proof}. 
Let $\alpha\notin \{0,1\}$, $t\neq 0$. There is a basis $\{e_1,\dots,e_n\}$ in $E$ so that 
$\langle e_i,e_j\rangle = \delta_{i,j}$. Let us show that $S(\alpha,t,E)$ is simple. 
If $I$ is a~nonzero ideal in $S(\alpha,t,E)$, then it contains a~nonzero element 
$x=\beta z_1 + \gamma z_2 + \sum\limits_{i=1}^n\alpha_i e_i$. 

If $\alpha_k\neq 0$ for some $k>0$, then $I$ contains an element 
$$
xe_k
= \alpha_k (z_1 + tz_2) + (\alpha\beta + (1-\alpha)\gamma ) e_k.
$$ 
Hence, $y=z_1+tz_2+\delta e_k\in I$ for $\delta = (\alpha\beta-\alpha\gamma+\gamma)/\alpha_k$. 
Then $I$ contains 
$$
(1-\alpha)yz_1-\alpha yz_2=(1-\alpha)z_1-\alpha tz_2.
$$
It means that $(1-\alpha)z_1\in I$ and $\alpha tz_2\in I$. Thus, $z_1,z_2\in I$. We also have
$$
I\ni z_1e_i = \alpha e_i
$$
for all $1\le i\le n$. Therefore, $I = S(\alpha,t,E)$.

If $\alpha_i=0$ for all $1\le i\le n$, then $0\neq \beta z_1 + \gamma z_2\in I$. 
It means that $z_i\in I$ for some $i\in\{1,2\}$. Then $I\ni z_ie_k$ and we have $e_k\in I$ 
for $1\le k\le n$ by assumptions. But it means that $e_k^2 = z_1+tz_2\in I$ and $z_1,z_2\in I$ as above. 

So, $I=S(\alpha,t,E)$ and $S(\alpha,t,E)$ is a simple algebra.

If $\alpha = 0$, then $Fz_1$ is a proper ideal in $S(\alpha,t,E)$. 
If $\alpha = 1$, then $Fz_2$ is a proper ideal in $S(\alpha,t,E)$. 
If $t=0$, then $Fz_1+\sum\limits_{i=1}^n Fe_i$ is a proper ideal in $S(\alpha,t,E)$. 
\hfill $\square$

We will use a notation $O(E)$ for a subgroup of $\Aut(S(\alpha,t,E))$ obtained by an extension 
of the orthogonal group of $E$, which elements fix $z_1$ and $z_2$.

{\bf Proposition 2}.
Let $E$ has a finite dimension $n\ge 2$ and $\alpha\notin \{0,1\}$, $t\neq 0$. 
\begin{itemize}
\item If $\alpha\neq 1/2$ or $t\neq \pm 1$, then $\Aut(S(\alpha,t,E)) \cong O(E)$.
\item If $\alpha = 1/2$ and $t=\pm 1$, then $\Aut(S(\alpha,t,E)) \cong \mathbb{Z}_2\times O(E)$.
\end{itemize}

{\sc Proof}. 
Let $A=S(\alpha,t,E)$ and $\varphi\in\mathrm{Aut}(A)$. There is a basis $\{e_1,\dots,e_n\}$ in $E$ 
so that $\langle e_i,e_j\rangle = \delta_{i,j}$. We want to describe all $x\in A$ with $\dim\mathrm{Ann}(x)=n$. 
Let $x\in A$ so that $\dim\mathrm{Ann}(x)=n$
and $x=\beta z_1 + \gamma z_2 + \alpha_1 e_1 +\ldots +\alpha_n e_n$. Then 
$$
xe_i = \alpha_i(z_1+tz_2) + (\beta\alpha + \gamma (1-\alpha)) e_i. 
$$
If $\beta\alpha+\gamma (1-\alpha) \neq 0$, then $xe_1,\dots, xe_n$ are linearly independent. Moreover, we have 
\begin{gather*}
xz_1 = \beta z_1 + \alpha (\alpha_1e_1+\dots + \alpha_n e_n), \\
xz_2 = \gamma z_2 + (1-\alpha) (\alpha_1e_1+\dots+\alpha_n e_n).
\end{gather*}
An assumption $\beta\alpha + \gamma(1-\alpha)\neq 0$ means that $\beta\neq 0$ or $\gamma\neq 0$. 
If $\beta\neq 0$, then $xz_1,xe_1,\dots,xe_n$ are linearly independent and $\dim\mathrm{Ann}(x)\le 1$, a contradiction. 
If $\gamma\neq 0$, then $xz_2,xe_1,\dots,xe_n$ are linearly independent and $\dim\mathrm{Ann}(x)\le 1$, a contradiction. 

So, $\beta\alpha + \gamma (1-\alpha) = 0$ and $\gamma = \alpha\beta/(\alpha - 1)$. 

Suppose that $\beta\neq 0$.
If $\alpha_i\neq 0$ for some~$i$, then $xe_i$, $xz_1$ and $xz_2$ are linearly independent 
and $\dim\mathrm{Ann}(x) < n$, a contradiction. 
Thus, $x=\beta (z_1 - \frac{\alpha}{1-\alpha}z_2)$. 
If $\beta = 0$, then $\gamma = 0$ and $x\in E$. 

Let us denote $U=F\cdot (z_1-\frac{\alpha}{1-\alpha}z_2)$. 
So, $x\in A$ and $\dim\mathrm{Ann}(x)=n$ if and only if $x\in E\cup U$ and $x\neq 0$. 
It means that $\varphi(x)\in E\cup U$ for any $x\in E\cup U$. 

Suppose that $\varphi(e)\in U$ for some $0\neq e\in E$. 
Since $\dim E \ge 2$, there exists $0\neq f\in E$ such that $\varphi (f)\in E$. 
It means that $\varphi(e+f)\notin E\cup U$, a contradiction. 
Hence, $\varphi (E)=E$, $\varphi(U)=U$.
Therefore,  
$\varphi(z_1-\frac{\alpha}{1-\alpha}z_2)=\delta(z_1-\frac{\alpha}{1-\alpha}z_2)$ 
for some nonzero~$\delta$.

We have $\varphi(z_1+z_2)=z_1+z_2$, so 
$$
\varphi(z_1)=(\alpha + \delta(1-\alpha))z_1+(\alpha-\delta\alpha)z_2.
$$
Since $z_1$ is an idempotent, hence $\varphi(z_1)$ is an idempotent. 
It means that 
$$
\alpha+\delta(1-\alpha),\alpha-\delta\alpha\in\{0,1\}.
$$
If $\alpha-\delta\alpha = 1$, then $\delta = (\alpha-1)/\alpha$. We have two cases:
\begin{enumerate}
\item $\alpha + \delta(1-\alpha) = 0$. 
It means that $\delta = \alpha/(\alpha - 1)$ and $\alpha = 1/2$ by above. 
Hence, $\delta = -1$ and $\varphi(z_1-z_2)=-z_1+z_2$. 
It is easy to see that $\varphi(z_1)=z_2$, $\varphi(z_2)=z_1$. We have 
$$
tz_1+z_2 = \varphi (z_1+tz_2) = \varphi(e_1^2) = \varphi(e_1)^2 = \gamma (z_1+tz_2)
$$
for some $\gamma\in F$. Hence, $\gamma = t$ and $t^2=1$, $t=\pm 1$. 
\item $\alpha + \delta(1-\alpha) = 1$. 
It means that $(\alpha - 1) = (\alpha - 1)\delta$ and $\delta = 1$.
\end{enumerate}

If $\delta = 1$, then 
$\varphi(z_1)=z_1$ and $\varphi(z_2)=z_2$. 
If $\delta = -1$, $\alpha = 1/2$, and $t=\pm 1$, then $\varphi(z_1)=z_2$ and $\varphi(z_2)=z_1$.
We have proved that $\varphi(e)\in E$ for any $e\in E$. 
In cases $\delta=1$ and $\delta=-1$, $\alpha=1/2$, $t=1$ the basis 
$e_1,\dots,e_n$ is mapped to an orthogonal basis of $E$, since $\varphi(z_1+tz_2)=z_1+tz_2$. 
It remains to prove that $\Aut(A) \cong \mathbb{Z}_2\times O(E)$, when
$\delta=t=-1$ and $\alpha=1/2$.
In the case, we have $\Aut(A) = \{ (\sigma,\psi) \mid \sigma\in S_2,\,
\langle \psi(e),\psi(f)\rangle = \mathrm{sgn}(\sigma) \langle e,f\rangle,\,e,f\in E \}$, 
where $\mathrm{sgn}(\sigma)$ denotes the sign of a~permutation $\sigma\in S_2$. 
Thus, $\pi\colon \Aut(A) \to \mathbb{Z}_2\times O(E)$ defined as follows,
$\pi( (\sigma,\psi) ) = (\sigma,\sqrt{-1}^{\mathrm{sgn}(\sigma)}\psi )$ is an isomorphism.
\hfill $\square$

\section{Sharped cubic form and associated Jordan algebra}

In this paragraph, we recall the results concerned sharped cubic forms and induced algebras, 
which occur to be Jordan. We follow the monograph~\cite{McCrimmon}.

A map $N\colon V\to F$ on a space $V$ is called cubic form if for any $\lambda\in F$ and $x,y\in V$
\[
N(x + \lambda y) = N(x) + \lambda N(x,y) + \lambda^2 N(y,x) + \lambda^3 N(x,y,z),
\]
where $N(x,y)$ is quadratic in $x$ and linear in $y$ and $N(x,y,z)$ is trilinear and symmetric.

Given a cubic form~$N$ on a space~$V$,
one can linearize it completely as follows,
$$
N(v,u,w) = N(v+u+w)-N(v+u)-N(v+w)-N(u+w)+N(u)+N(v)+N(w).
$$
Hence, $N(r,r,r) = 6N(r)$.

{\bf Definition 1}.
Let $V$ be a vector space endowed with a cubic form~$N$
and let $c\in V$ be such that $N(c) = 1$ (we call $c$ as basepoint).
Denote by $N(x,y,z)$ the complete linearization of~$N$.
We introduce quadratic spur function, linear trace form, two bilinear forms 
and give another definition of the form $N(x,y)$:
\begin{gather}
S(r) = N(r,r,c)/2, \quad
T(r) = N(r,c,c)/2, \label{def:S,T} \\
S(r,q) = N(r,q,c), \quad N(r,q) = N(r,r,q)/2, \label{def:S2,N2} \\
(r,q) = T(r)T(q) - S(r,q). \label{def:(,)J}
\end{gather}
From the definition we derive that
\begin{equation} \label{T,SOnUnit}
S(c) = T(c) = 3, \quad
(r,c) = T(r).
\end{equation}

Another way to define the function $N(r,q)$ from the given norm function~$N$ is to present
\begin{equation} \label{Norm->Norm2}
N(r+tq) = N(r) + tN(r,q) + t^2N(q,r) + t^3N(q).
\end{equation}
Thus, $N(r,q,s)$ is a linearization of $N(x,y)$.

{\bf Definition 2}.
Let $V$ be a vector space endowed with a cubic form~$N$
and basepoint~$c$. A sharp map~$\#$ on $V$ for $(N,c)$ is a quadratic operator on~$V$ 
satisfying the following relations:
\begin{gather}
(r^{\#},q) = N(r,q), \\
(r^{\#})^{\#} = N(r)r, \\
c_{\#}r = T(r)c - r,
\end{gather}
where the sharp-product is given by the formula
\begin{equation} \label{sharpproduct}
r_\# q = (r+q)^\# - r^\# - q^\#.
\end{equation}
Under the conditions, we call $(N,\#,c)$ as a sharped cubic form.

Due to the definitions, $c^{\#} = c$.

{\bf Theorem 1}.
Let $V$ be a vector space endowed with a sharped cubic form $(N,\#,c)$. Then

a) $V$ under the product
\begin{equation} \label{product}
rq = \frac{1}{2}(r_\# q + T(r)q + T(q)r - S(r,q)c)
\end{equation}
is an algebra with a unit~$c$.
Moreover, every $r\in V$ satisfies the cubic identity
\begin{equation} \label{cubic-identity}
r^3 - T(r)r^2 + S(r)r - N(r)c = 0,
\end{equation}
and $r^{\#} = r^2 - T(r)r + S(r)c$ (so, $r^{\#}r = N(r)c$).

b) $V$ is Jordan and the operators
\begin{gather}
U_r(s) = (r,s)r - {r^\#}_\# s, \quad
U_{r,q}(s) := U_{r+q}(s) - U_r(s) - U_q(s), \label{def:Uoperator} \\
\{r,s,q\} := (r,s)q + (q,s)r - (r_\# q)_\# s \label{tripleproduct}
\end{gather}
coincide with the classical $U$-operators and Jordan triple product respectively.

In~\cite{McCrimmon}, Theorem 1b) is proved via the following formulas.

{\bf Proposition 3}.
Let $V$ be a vector space endowed with a sharped cubic form $(N,\#,c)$. 
Then the following identities hold on~$V$:
\begin{gather}
S(r) = T(r^{\#}), \quad
S(r,q) = T(r_{\#}q), \\
(rq,s) = (r,qs), \quad
(r_{\#}q,s) = (r,q_{\#}s) = N(r,q,s), \quad
(U_r(s),q) = (s,U_r(q)), \allowdisplaybreaks \\
r^{\#}{}_{\#}(q_{\#}r) = N(r)q + (r^{\#},q)r, \quad
(r^{\#}{}_{\#}q)_{\#}r = N(r)q + (r,q)r^{\#}, \\
(r_{\#}q)^{\#} + r^{\#}{}_{\#}q^{\#} = (r^{\#},q)q + (r,q^{\#})r, \\
r^{\#}{}_{\#}r = -T(r)r^{\#} - T(r^{\#})r + (S(r)T(r)-N(r))c, \\
S(r^{\#},r) = S(r)T(r) - 3N(r), \quad
(r^{\#},r) = 3N(r), \\
U_r(c) = r^2, \quad
U_{r,q}(c) = 2rq, \quad
U_r(r^{\#}) = N(r)r, \quad
(U_r(q))^{\#} = U_{r^{\#}}(q^{\#}), \\
U_r U_{r^{\#}} = N(r)^2 \mathrm{id}, \quad
\{r,r^{\#},q\} = 2N(r)q, \\
N(U_r(q)) = N(r)^2N(q), \quad N(r^{\#}) = N(r)^2.
\end{gather}

\section{Generalized sharped cubic form and its algebra}

Now, we suggest a construction, which generalizes sharped cubic form.

{\bf Definition 3}.
Let $V$ be a vector space endowed with a cubic form~$N$
and let $c\in V$ be such that $N(c) = 1$.
We also assume that a symmetric bilinear form~$\Delta$ is defined on~$V$ in such manner that 
\begin{equation} \label{delta(r,c)=0}
\Delta(r,c) = 0 
\end{equation}
for all~$r\in V$.
Denote by $N(x,y,z)$ the complete linearization of~$N$.
As in Definition~1, we introduce quadratic spur function $S(r)$, linear trace form~$T(r)$, 
bilinear form~$S(r,q)$ by~\eqref{def:S,T}, quadratic in $r$ and bilinear in $q$ 
form~$N(r,q)$ by~\eqref{def:S2,N2} and new bilinear form
\begin{equation} \label{GNorm}
(r,q) = T(r)T(q) - S(r,q) - \Delta(r,q).
\end{equation}
Let us call a pair~$(N,\Delta)$ as a generalized cubic form.
When $\Delta = 0$, we have an ordinary cubic form.

The equalities~\eqref{T,SOnUnit} follow from the definition immediately.

{\bf Definition 4}.
Let $V$ be a vector space endowed with a generalized cubic form~$(N,\Delta)$ and basepoint~$c$. 
A sharp map~$\#$ on $V$ for $(N,\Delta,c)$ is a quadratic operator on~$V$ 
satisfying the following relations:
\begin{gather}
(r_\#q,r) + (r^\#,q) = 3N(r,q), \label{def:GSharp1} \\
(r^\#)^\#
= (N(r) + \Delta(r^\#,r))r, \label{def:GSharp2} \\
c_{\#}r = T(r)c - r, \label{def:GSharp3}
\end{gather}
where the sharp-product is given by~\eqref{sharpproduct}.
Under the conditions, we call $(N,\Delta,\#,c)$ as a~generalized sharped cubic form.

As above, we conclude that $c^{\#} = c$.

Given a~generalized sharped cubic form~$(N,\Delta,\#,c)$, let us
define a product on~$V$ by~\eqref{product}.

{\bf Proposition 4}.
Fix a scalar $\lambda\in F\setminus\{-1\}$.
Given a cubic form $N$ defined on a vector space~$V$ with a basepoint~$c$, 
we get a generalized cubic form by the formulas
$$
(r,q)
= \frac{1}{\lambda+1}\left((1+\lambda/3)T(r)T(q)-S(r,q)\right), \quad 
\Delta(r,q) = \lambda \left((r,q) - \frac{T(r)T(q)}{3}\right).
$$

{\sc Proof}.
The identity~\eqref{GNorm} holds by the definition.
Also, we get $(r,c) = \frac{T(r)+\lambda T(r)}{1+\lambda} = T(r)$, hence, $\Delta(r,c) = 0$.
\hfill $\square$

Let us call a generalized cubic form defined in Proposition~4 as an inner one.

{\bf Example 1}. 
Let $A$ be an associative commutative algebra over a field~$F$ generated by the unit of $F$ 
and an element $\lambda$ such that $\lambda^2 = 0$.
Consider the space $V = A\otimes_F F^3\cong A^{\otimes3}$. 
We endow~$V$ with the form
$N((x,y,z)) = xyz$ and take $c = (1,1,1)$. 
Formally, $N$ is not a~cubic form, since it maps $A^{\otimes3}$ to $A$ instead of $F$. 
However, we get an interesting example of an algebra.
Denote $r = (x,y,z)$, $q = (x',y',z')$. Then 
\begin{gather*}
T(r) = x+y+z, \quad
S(r) = xy + xz + yz, \\
S(r,q) = x(y'+z') + y(x'+z') + z(x'+y'), \
N(r,q) = xyz' + xy'z + x'yz.
\end{gather*}
We define 
$$
r^\# = (yz,xz,xy) - \lambda(y^2+z^2+2x(y+z),x^2+z^2+2y(x+z),x^2+y^2+2z(x+y)).
$$
and get by~\eqref{product},
$$
rq
= (1+\lambda)(xx',yy',zz') + \lambda(yz'+y'z,xz'+x'z,xy'+x'y) - \lambda T(r)T(q)c.
$$
By the definition, 
$$
c^\# = (1-6\lambda)c, \quad
rc = r - 2\lambda T(r)c,\quad 
c_{\#}r = (1-4\lambda)T(r)c - r.
$$

Now, we consider two pairs of bilinear maps $(\cdot,\cdot)$ and $\Delta$ such that~\eqref{GNorm} holds:
$$
(r,q) = xx'+yy'+zz' + 3\lambda S(r,q), \quad 
\Delta(r,q) 
= - 3\lambda S(r,q).
$$
When $\Delta = \lambda = 0$, we get an ordinary sharped cubic form.

\newpage

It is easy to check that a cubic form $(N,\Delta,\#,c)$ satisfies~\eqref{def:GSharp1}. 
Instead of~\eqref{def:GSharp2}, the following relation holds:
$$
(r^\#)^\# = (N(r) + \Delta(r^\#,r) + 2\lambda T(r)S(r))r
+ 2\lambda S(r)r^\# - 2\lambda (T(r)N(r)+S(r)^2)c. 
$$

There is some routine (see the code in GAP~\cite{Code}) to derive the relations fulfilled 
for such maps (see the same or close equalities derived for all generalized sharped cubic forms below):
\begin{gather*}
\Delta(r,c) = -6\lambda T(r), \quad
(r,c) = (1+6\lambda)T(r), \quad 
(r^\#,r) = 3N(r), \\
T(r^\#) = (1-6\lambda)(S(r) - \Delta(r,r)) - 2\lambda T(r)^2, 
\allowdisplaybreaks \\
(r_\#q,s) - (r,q_\#s) 
= \Delta(r,q_\#s)  - \Delta(r_\#q,s)
= T(r)\Delta(q,s) - T(s)\Delta(r,q), \\ 
(r_\#q,s) + (q_\#s,r) + (s_\#r,q)  = 3N(r,q,s), \\
N(r^\#) = N(r)(N(r) + \Delta(r^\#,r)).
\end{gather*}

{\bf Example 2}.
The general case of the split spin factor $S(\alpha,t,E)$ with 
\begin{gather}
N(az_1 {+} bz_2 {+}v) = ab(\alpha a+\bar{\alpha}b)
- \langle v,v\rangle(\bar{\alpha}ta+\alpha b), \\
\Delta(az_1 + bz_2 + v,kz_1 + lz_2 + u) 
= \alpha(\alpha-1)(a-b)(k-l)-\langle u,v\rangle(\bar{\alpha}+\alpha t),\\
(az_1 + bz_2 + v)^\# 
= (\alpha a+\bar{\alpha}b)(bz_1+az_2) 
+ (t-1)\langle v,v\rangle(-\bar{\alpha}z_1+\alpha z_2) 
- (\bar{\alpha}a {+} \alpha b)v, \label{SplitSpinSharp}
\end{gather} 
where $\bar{\alpha} = 1-\alpha$ and $c = z_1 + z_2$
is a generalized sharped cubic form. 
Here $a,b,k,l\in F$ and $u,v\in E$.

Indeed, $N(c) = 1$, $\Delta(r,c) = 0$. Further, for $r = az_1 + bz_2 + v$, we have
\begin{multline*}
2T(r)
= N(r+2c) - 2N(r+c) - N(2c) + 2N(c) + N(r) \\
= (a+2)(b+2)(\alpha(a+2)+\bar{\alpha}(b+2)) - 2(a+1)(b+1)(\alpha(a+1)+\bar{\alpha}(b+1)) \\
+ ab(\alpha a+\bar{\alpha}b) - 6 
= 2( (1+\alpha)a + (2-\alpha)b ),
\end{multline*}
and~\eqref{def:GSharp3} holds, since
\begin{multline*}
r_\# c 
= (r+c)^\# - r^\# - c^\# \\
= (\alpha(a+1)+\bar{\alpha}(b+1))((b+1)z_1+(a+1)z_2)
- (\alpha a+\bar{\alpha}b)(bz_1+az_2) - v - c \\
= ( (1+\alpha)a + (2-\alpha)b )(z_1+z_2) - az_1 - bz_2 - v
= T(r)c - r.
\end{multline*}

Now, we compute due to the definition,
$$
\Delta(r,r^\#)
= \alpha(\alpha-1)(a-b)( (\alpha a+\bar{\alpha}b)(b-a) - (t-1)\langle v,v\rangle ) 
+ (\bar{\alpha}a +\alpha b)\langle v,v\rangle(\bar{\alpha}+\alpha t).
$$

It is not difficult to show that
$$
N(r) + \Delta(r,r^\#)
= (\bar{\alpha}a+\alpha b)( (\alpha a+\bar{\alpha}b)^2
+ (2\alpha-1)(t-1)\langle v,v\rangle).
$$
On the other hand, the projection of $(r^\#)^\#$ on $E$ equals by~\eqref{SplitSpinSharp}
$$
(\bar{\alpha}a + \alpha b)\big( (\alpha a+\bar{\alpha}b)(\bar{\alpha}b + \alpha a) 
+ (2\alpha-1)(t-1)\langle v,v\rangle \big)v 
= (N(r) + \Delta(r,r^\#))v.
$$
We leave the check that the coordinates of $(r^\#)^\#$ at $z_1$ and $z_2$ 
also equal to the ones of $(N(r) + \Delta(r,r^\#))r$.
Thus, \eqref{def:GSharp2} follows.

Finally, we have to derive~\eqref{def:GSharp1}. 
For this, we write down the required forms for $r = az_1 + bz_2 + v$ and $s = kz_1 + lz_2 + u$:
\begin{multline*}
N(r,s)
= (N(2r+s) - 2N(r+s) - N(2r) + 2N(r) + N(s))/2 \\
=  \alpha a^2 l + 2\alpha ab k + 2\bar{\alpha}abl + \bar{\alpha}b^2k 
- \langle v,v\rangle(\bar{\alpha}t k+\alpha l)
- 2\langle v,u\rangle(\bar{\alpha}t a+\alpha b),
\end{multline*}

\vspace{-0.75cm}

\begin{multline*}
S(r,s)
= N(r+s+c) - N(r+s) - N(r+c) - N(s+c) + N(r) + N(s) + N(c) \\
= 2(\alpha ak+al+bk+\bar{\alpha}bl) -2(\bar{\alpha}t+\alpha)\langle v,u\rangle,
\end{multline*}

\vspace{-0.75cm}

\begin{equation} \label{inner-product-derived}
(r,s) 
= T(r)T(s) - S(r,s) - \Delta(r,s)
= (1+\alpha)ak + (2-\alpha)bl + (1+\alpha+(2-\alpha)t)\langle v,u\rangle,
\end{equation}

\vspace{-0.5cm}

\begin{multline*}
r_\# s
= (r+s)^\# - r^\# - s^\#
= (\alpha a + \bar{\alpha}b)(lz_1+kz_2)
+ (\alpha k + \bar{\alpha}l)(bz_1+az_2) \\
+ 2(t-1)(-\bar{\alpha}z_1+\alpha z_2)\langle v,u\rangle
- (\bar{\alpha}a+\alpha b)u
- (\bar{\alpha}k+\alpha l)v.
\end{multline*}
Applying these formulas to compute $(r_\#s,r)$ and $(r^\#,s)$, 
we prove~\eqref{def:GSharp1}.

We write down the product on $S(\alpha,t,E)$ by~\eqref{product}:
\begin{multline} \label{product-derived}
rs
= \frac{1}{2}(r_\#s + T(r)s + T(s)r - S(r,s)c)
= akz_1+blz_2+\langle v,u\rangle(z_1 + tz_2) \\
+ (\alpha k +\bar{\alpha}l)v + (\alpha a+\bar{\alpha}b)u,
\end{multline}
which, up to rescalling the bilinear form on $E$, coincides with the initial product on $S(\alpha,t,E)$.

{\bf Remark 1}.
Denote $\lambda = \dfrac{3\alpha(1-\alpha)}{(1+\alpha)(\alpha-2)}$. 
Then the generalized sharped cubic form on the split spin factor $S(\alpha,E)$ 
is inner with $\lambda$, since for $r = az_1+bz_2+v$ and $s = kz_1 + lz_2 + u$ we have
\begin{multline*}
(r,s) - T(r)T(s)/3 \\
= (az_1 + bz_2 + v,kz_1 + lz_2 + u) - \frac{( (1+\alpha)a 
	+ (2-\alpha)b )( (1+\alpha)k + (2-\alpha)l )}{3} \\
= \frac{1}{3} \bigg( ak(1+\alpha)(3-(1+\alpha)) + bl(2-\alpha)(3-(2-\alpha)) 
 - (1+\alpha)(2-\alpha)(al+bk) + \frac{3(1+\alpha)}{\alpha}\langle v,u\rangle \bigg) \\
= \frac{(1+\alpha)(2-\alpha)}{3}\left( (a-b)(k-l) + \frac{3\langle v,u\rangle}{\alpha(2-\alpha)} \right)
= \frac{1}{\lambda} \Delta(r,s).
\end{multline*}

Let us return to generalized sharped cubic forms and relations concerned with them.

{\bf Lemma 1}.
Given a vector space~$V$, let $(N,\Delta,\#,c)$ be a~generalized sharped cubic form on~$V$. 
Then the following identities hold:
\begin{gather}
T(r^\#) = S(r) - \Delta(r,r), \label{T(sharp(r))S(r)} \\
(r^\#,r) = 3N(r), \label{(sharp(r),r)} \\
S(r^\#,r) = T(r)(S(r)-\Delta(r,r)) - 3N(r) - \Delta(r^\#,r), \label{S(sharp(r),r)} \\
{r^\#}_\#(r_\# q) = (N(r)+\Delta(r^\#,r))q + (N(r,q)+\Delta(r^\#,q)+\Delta(r,r_\#q))r, \label{Adjoint'}
\end{gather}

\vspace{-1.0cm}

\begin{multline} \label{Adjoint''}
(r_{\#}q)^\#+{r^\#}_{\#}q^\#
= (N(q,r)+\Delta(q,r_\#q)+\Delta(r,q^\#))r \\
+ (N(r,q)+\Delta(r^\#,q)+\Delta(r,r_\#q))q,
\end{multline}
\begin{equation} \label{sharpproduct(sharp(r),r)}
{r^\#}_\# r
= -T(r)r^\# - T(r^\#)r + (T(r)(S(r)-\Delta(r,r))-N(r) - \Delta(r^\#,r))c.
\end{equation}

{\sc Proof}.
Involving~\eqref{GNorm},~\eqref{def:GSharp1},~\eqref{def:GSharp3}, and~\eqref{T,SOnUnit}, we get
$$
T(r^\#)
= (r^\#,c)
= 3N(r,c) - (r_\#c,r)
= 3S(r) - T(r)^2 + (r,r)
= S(r) - \Delta(r,r).
$$

By~\eqref{def:GSharp1}, we rewrite
$$
9N(r) 
= 3N(r,r)
= (r_\#r,r) + (r^\#,r)
= 3(r^\#,r),
$$
since $r_{\#}r = (2r)^{\#} - r^\# - r^\# = 2r^\#$.
Thus, we derive~\eqref{(sharp(r),r)}.

Applying~\eqref{GNorm} and already proved formulas, we deduce
$$
S(r^\#,r)
= T(r)T(r^\#) - (r^\#,r) - \Delta(r^\#,r)
= T(r)(S(r)-\Delta(r,r)) - 3N(r) - \Delta(r^\#,r).
$$

Let us put $r+tq$ instead of~$r$ into~\eqref{def:GSharp2}, where $t\in F$.
Joint with~\eqref{Norm->Norm2}, we get
\begin{multline*}
((r+tq)^{\#})^{\#}
= (r^{\#}+t^2q^{\#}+tr_{\#}q)^{\#}
= (r^{\#}+tr_{\#}q)^{\#} + t^4(q^{\#})^{\#}
+ t^2(r^{\#}+tr_{\#}q)_{\#}q^{\#} \\
= (r^{\#})^{\#} + t^2(r_\#q)^{\#}+t{r^{\#}}_{\#}(r_{\#}q)+t^4(q^{\#})^{\#}
 +t^2{r^\#}_{\#}q^\#+t^3(r_{\#}q)_\#q^\#;
\end{multline*}

\vspace{-1.0cm}

\begin{multline*}
(N(r+tq)+\Delta(r+tq,(r+tq)^\#))(r+tq) \\
= (N(r) + tN(r,q) + t^2N(q,r) + t^3N(q))(r+tq)
+ \Delta(r+tq,r^\#+t^2q^\#+tr_\#q)(r+tq).
\end{multline*}
Comparing coefficients at~$t$ and at~$t^2$, we 
derive~\eqref{Adjoint'} and~\eqref{Adjoint''} respectively.

Finally, we apply~\eqref{def:GSharp3} twice and 
then~\eqref{T(sharp(r))S(r)},~\eqref{Adjoint'} with $q=c$:
\begin{multline*}
{r^\#}_\# r
= {r^\#}_\# (T(r)c - r_\# c)
= T(r)(T(r^\#)c-r^\#) - {r^\#}_\# (r_\# c) \\
= T(r)(S(r)-\Delta(r,r))c - T(r)r^\#
- (N(r)+\Delta(r^\#,r))c - (N(r,c)+\Delta(r,T(r)c-r))r \\
= (T(r)(S(r)-\Delta(r,r)) - N(r) - \Delta(r^\#,r) )c
- T(r)r^\# - T(r^\#)r. \qquad \square
\end{multline*}

{\bf Theorem 2}.
Given a vector space~$V$, let $(N,\Delta,\#,c)$ be a~generalized sharped cubic form on~$V$.
Define a product~$\cdot$ on~$V$ by~\eqref{product}.
Then the conclusion of Theorem 1a) holds for $(V,\cdot)$.

{\sc Proof}.
Let us check that $c$ is a unit 
by~\eqref{def:S,T},~\eqref{def:S2,N2},~\eqref{T,SOnUnit},~\eqref{def:GSharp3}:
$$
2rc = r_{\#}c + T(r)c + T(c)r - S(r,c)c
= 2T(r)c - r + 3r - 2T(r)c = 2r.
$$

Due to the definitions,
$$
r^2 = \frac{1}{2}(r_{\#}r + 2T(r)r - S(r,r)c)
= r^\# + T(r)r - S(r)c.
$$

It remains to show~\eqref{cubic-identity}.
For this, we apply~\eqref{S(sharp(r),r)},~\eqref{sharpproduct(sharp(r),r)}:
\begin{multline*}
2rr^\#
= r_\# r^\# + T(r)r^\# + T(r^\#)r - S(r,r^\#)c \\
= (T(r)(S(r)-\Delta(r,r)) - N(r) - \Delta(r^\#,r) )c
- (T(r)(S(r)-\Delta(r,r)) - 3N(r) - \Delta(r^\#,r))c \\
= 2N(r)c,
\end{multline*}
i.\,e. $rr^\# = N(r)c$, which is equivalent to~\eqref{cubic-identity}.
\hfill $\square$

Let us state further relations fulfilled on an algebra endowed with a generalized sharped cubic.

{\bf Lemma 2}.
Given a vector space~$V$, let $(N,\Delta,\#,c)$ be a~generalized sharped cubic form on~$V$. 
Then the following identities hold:
\begin{gather}
T(r_\#q) = S(r,q) - 2\Delta(r,q), \label{T(sharpproduct(r,q))} \\
N(r^\#) = N(r)(N(r) + \Delta(r^\#,r)), \label{Norm(sharp(r))} \\
(r_\#q,s) + (q_\#s,r) + (s_\#r,q)  = 3N(r,q,s), \label{InvFormUnderSharpCycle} \\
{r^\#}_\# r^\# = 2(N(r)+\Delta(r^\#,r)))r, \label{sharpproduct(sharp(r),sharp(r))}\\
U_c(r) = r, \quad
U_r(c) = r^2 + \Delta(r,r)c, \quad
\frac{1}{2}U_{r,q}(c) = rq + \Delta(r,q)c, \label{UopIdentity} \\
U_r(r) = r^3 - 2\Delta(r,r)r + (\Delta(r^\#,r) + T(r)\Delta(r,r))c, \label{U_r(r)}  \\
U_r(r^\#) = (N(r)-2\Delta(r^\#,r))r. \label{U_r(sharp(r))}
\end{gather}

{\sc Proof}.
By~\eqref{T(sharp(r))S(r)} and the definition of the sharp product, 
we get~\eqref{T(sharpproduct(r,q))}:
\begin{multline*}
T(r_\#q)
= T((r+q)^\#) - T(r^\#) - T(q^\#) \\
= S(r+q) - \Delta(r+q,r+q) - S(r) + \Delta(r,r) - S(q) + \Delta(q,q)
= S(r,q) - 2\Delta(r,q).
\end{multline*}

We use~\eqref{(sharp(r),r)} and the axiom~\eqref{def:GSharp2} to compute
$$
3N(r^\#)
= ((r^\#)^\#,r^\#)
= ((N(r) + \Delta(r^\#,r))r,r^\#)
= 3(N(r) + \Delta(r^\#,r))N(r),
$$
hence,~\eqref{Norm(sharp(r))} is proved.

Linearization of~\eqref{def:GSharp1} implies~\eqref{InvFormUnderSharpCycle}.
Since $q_\# q = 2q^\#$, the equality~\eqref{sharpproduct(sharp(r),sharp(r))} 
follows from~\eqref{def:GSharp2}.

By the definition of the $U$-operator,~\eqref{T,SOnUnit},~\eqref{def:GSharp3}, 
and~\eqref{T(sharp(r))S(r)}, we have
\begin{gather*}
U_c(r) 
= (c,r)c - {c^\#}_\# r 
= T(r)c - c_\# r
= r, \\
U_r(c)
= (c,r)r - {r^\#}_\# c
= T(r)r - T(r^\#)c + r^\#
= r^2 - T(r^\#)c + S(r)c
= r^2 + \Delta(r,r)c,
\end{gather*}
and the third equality from~\eqref{UopIdentity} follows immediately.

Applying~\eqref{cubic-identity},~\eqref{GNorm},~\eqref{T(sharp(r))S(r)},~\eqref{sharpproduct(sharp(r),r)}, 
we derive~\eqref{U_r(r)}:
\begin{multline*}
U_r(r) 
- r^3 + 2\Delta(r,r)r - (\Delta(r^\#,r) + T(r)\Delta(r,r))c \\
= (r,r)r - {r^\#}_\# r - T(r)r^2 + S(r)r - N(r)c + 2\Delta(r,r)r - (\Delta(r^\#,r) +  T(r)\Delta(r,r))c \\
= (r,r)r + T(r)r^\# + T(r^\#)r - (T(r)(S(r)-\Delta(r,r))-N(r) - \Delta(r^\#,r))c \\
- T(r)r^2 + S(r)r - N(r)c + 2\Delta(r,r)r - (\Delta(r^\#,r) +  T(r)\Delta(r,r))c \\
= ((r,r) + S(r) - T(r)^2 + \Delta(r,r))r
= 0.
\end{multline*}

Finally, the formula~\eqref{U_r(sharp(r))} holds by~\eqref{(sharp(r),r)} 
and~\eqref{sharpproduct(sharp(r),sharp(r))}.
\hfill $\square$

\newpage

{\bf Corollary 1}.
Given a vector space~$V$, let $(N,\Delta,\#,c)$ be a~generalized sharped cubic form on~$V$.
Then $(r,s) = T(rs)$ for all $r,s\in V$.

{\sc Proof}.
By the definition~\eqref{product}, we write down
$$
2T(rs)
= T(r_\#s) + 2T(r)T(s) - S(r,s)T(c)
\mathop{=}\limits^{\eqref{T(sharpproduct(r,q))}}
2(T(r)T(s) - S(r,s) - \Delta(r,s))
\mathop{=}\limits^{\eqref{GNorm}}
2(r,s).
$$
Corollary is proved.
\hfill $\square$

{\bf Lemma 3}.
Given a vector space~$V$, let $(N,\Delta,\#,c)$ be a~generalized sharped cubic form on~$V$.
Then the following identities are equivalent: 
\begin{gather}
(r_\#q,s) = N(r,q,s) + \frac{1}{3}(T(r)\Delta(q,s)+T(q)\Delta(r,s) - 2T(s)\Delta(r,q)), 
 \label{(sharpproduct(r,q),s)} \\
(r_\#q,s) = (r,q_\#s) + T(r)\Delta(q,s) - T(s)\Delta(r,q), \label{InvFormUnderSharpprod} \\
(rq,s) = (r,qs). \label{InvForm}    
\end{gather}

{\sc Proof}.
Note that~\eqref{(sharpproduct(r,q),s)} immediately implies~\eqref{InvFormUnderSharpprod} 
by~\eqref{InvFormUnderSharpCycle}.
Suppose that~\eqref{InvFormUnderSharpprod} holds. 
Then by~\eqref{InvFormUnderSharpCycle}, we have
\begin{multline*}
3N(r,q,s)
= (r_\#q,s) + (q_\#s,r) + (s_\#r,q) \\
= 3(r_\#q,s) - T(r)\Delta(q,s) - T(q)\Delta(r,s) + 2T(s)\Delta(r,q),
\end{multline*}
and~\eqref{(sharpproduct(r,q),s)} is fulfilled.

Now, we prove the last equivalency:
\begin{multline*}
2(rq,s) - 2(r,qs)    
= (r_\#q,s) + T(r)(q,s) + T(q)(r,s) - S(r,q)(c,s) \\
- (r,q_\#s) - T(q)(r,s) - T(s)(r,q) + S(q,s)(r,c) \\
= (r_\#q,s) - (r,q_\#s) - T(r)\Delta(s,q) + T(s)\Delta(r,q) \\
+ T(r)((q,s)+\Delta(q,s)+S(q,s)) 
- T(s)((r,q)+\Delta(r,q)+S(r,q)) \\
\mathop{=}\limits^{\eqref{GNorm}}
(r_\#q,s) - (r,q_\#s) - T(r)\Delta(s,q) + T(s)\Delta(r,q) + T(r)T(q)T(s)-T(s)T(r)T(q).
\end{multline*}
Hence,~\eqref{InvFormUnderSharpprod} and~\eqref{InvForm} are equivalent.
\hfill $\square$

{\bf Remark 2}.
Given a vector space~$V$, let $(N,\Delta,\#,c)$ be a generalized cubic form 
satisfying~\eqref{def:GSharp2} and~\eqref{def:GSharp3}.
Then $V$ is a generalized sharped cubic form with invariant form $(\cdot,\cdot)$ 
if and only if $V$ satisfies
\begin{equation}
(r^\#,q) = N(r,q) + (T(r)\Delta(r,q) - T(q)\Delta(r,r))/3, \label{def:GSharp1-old}
\end{equation}
for all $r,q\in V$.
Indeed, suppose that $V$ satisfies~\eqref{def:GSharp1-old}. 
Then a linearization of~\eqref{def:GSharp1-old} implies~\eqref{(sharpproduct(r,q),s)}.
Taking~\eqref{(sharpproduct(r,q),s)} with $s = r$, we get
$$
(r_\#q,r) = 2N(r,q) + \frac{1}{3}(T(q)\Delta(r,r)-T(r)\Delta(r,q)),
$$
which sum with~\eqref{def:GSharp1-old} gives~\eqref{def:GSharp1}.
Conversely, the identity~\eqref{def:GSharp1-old} holds on every generalized sharped cubic form 
with invariant form~$(\cdot,\cdot)$, since~\eqref{def:GSharp1-old} 
follows from~\eqref{(sharpproduct(r,q),s)} with $q=r$.

Now we state identities concerned the triple product~$\{r,s,q\}$ defined by~\eqref{tripleproduct}.

{\bf Lemma 4}.
Given a vector space~$V$, let $(N,\Delta,\#,c)$ be a~generalized sharped cubic form on~$V$.
Then the following identities hold:
\begin{equation}\label{triple(r,r,q)} 
\{r,r,q\} {=} (2r^2{-}\Delta(r,r))q {-} 3\Delta(r,q)r 
{+} \big(2T(r)\Delta(r,q) {-} (r^\#,q) {+} \Delta(r_\# q,r) {+} N(r,q)\big)c,
\end{equation}

\vspace{-0.75cm}

\begin{multline} \label{triple(r,s,q)+triple(s,r,q)}
\{r,s,q\} + \{s,r,q\} = (4(rs) - 2\Delta(r,s))q -3\Delta(r,q)s - 3\Delta(s,q)r \\
 + \big(2T(r)\Delta(s,q) + 2T(s)\Delta(r,q) - (r_\# s,q) + \Delta(s_\# q, r) 
 + \Delta(r_\# q, s) + N(r,s,q)\big)c,
\end{multline}

\vspace{-1cm}

\begin{multline} \label{(r,s,q)}
(r,s,q)
= \frac{1}{4}\big( \{s,r,q\} - \{s,q,r\}
+ \Delta(q,s)r - \Delta(r,s)q \\
-(\Delta(q_\#s,r) - \Delta(r_\#s,q) + 2T(r)\Delta(s,q) 
- 2T(q)\Delta(r,s) - (r_\# s,q) + (r,s_\# q))c\big).
\end{multline}

{\sc Proof}.
Denote 
$$
\eta(r,q) = (2T(r)\Delta(r,q) + T(q)\Delta(r,r) + \Delta(r^\#,q) 
 + \Delta(r_\# q,r) + N(r,q) - T(r)S(r,q))c.
$$
Then linearization of~\eqref{sharpproduct(sharp(r),r)} gives
\begin{multline*}
0 = (r_\# q)_\# r + {r^\#}_\# q + T(q)r^\#+T(r)(r_\# q) 
 + T(r^\#)q+T(r_\# q)r - T(q)S(r)c + \eta(r,q).
\end{multline*}
By \eqref{product},~\eqref{T(sharp(r))S(r)}, and~\eqref{T(sharpproduct(r,q))}, we have
\begin{multline*}
0 = (r_\# q)_\# r + r^2_\# q + S(r)(c_\# q) + T(q)(r^\# - S(r)c) + S(r)q - \Delta(r,r)q + S(r,q)r \\
- 2\Delta(r,q)r + \eta(r,q).
\end{multline*}
The identity~\eqref{def:GSharp3} implies
\begin{multline*}
0 = (r_\# q)_\# r + r^2_\# q + S(r)T(q)c + T(q)r^2 -T(r)T(q)r + S(r,q)r - 2\Delta(r,q)r \\ 
- \Delta(r,r)q + \eta(r,q) 
= (r_\# q)_\# r + r^2_\# q + S(r)T(q)c + T(q)r^2 - T(r)T(q)r \\
+ S(r,q)r - 2\Delta(r,q)r - \Delta(r,r)q + \eta(r,q).
\end{multline*}

The identity $T(r^2)=(r,r)$ and \eqref{product} imply
\begin{multline*}
0 = (r_\# q)_\# r + 2r^2q - T(r^2)q - T(q)r^2 + S(r^2,q)c + S(r)T(q)c \\ 
+ T(q)r^2 - T(r)T(q)r + S(r,q)r - 2\Delta(r,q)r - \Delta(r,r)q + \eta(r,q) \\
= (r_\# q)_\# r + 2r^2q + S(r^\#,q)c - (r,r)q + T(r)S(r,q)c - S(r)S(c,q)c \\
+ S(r)T(q)c - T(r)T(q)r + S(r,q)r - 2\Delta(r,q)r - \Delta(r,r)q + \eta(r,q). 
\end{multline*}
By \eqref{tripleproduct} and \eqref{GNorm}, we have
\begin{multline*}    
0 = -\{r,r,q\} + (r,q)r + 2r^2q - T(r)T(q)r +S(r,q)r - S(r)T(q)c + S(r^\#,q)c + T(r)S(r,q)c \\
 - 2\Delta(r,q)r - \Delta(r,r)q + \eta(r,q) 
 = -\{r,r,q\} + 2r^2q - 3\Delta(r,q)r - \Delta(r,r)q \\
 + \big(2T(r)\Delta(r,q) + T(q)\Delta(r,r) + \Delta(r^\#,q) 
 + \Delta(r_\# q,r) + N(r,q) - S(r)T(q) + S(r^\#,q)\big)c.
\end{multline*}
By~\eqref{GNorm}, we have
\begin{multline*}
0 = -\{r,r,q\} + (2r^2-\Delta(r,r))q - 3\Delta(r,q)r \\
+ \big(2T(r)\Delta(r,q) + T(q)\Delta(r,r)
+ T(r^\#)T(q) - (r^\#,q) + \Delta(r_\# q,r) + N(r,q) - S(r)T(q) \big)c.
\end{multline*}
The relation~\eqref{T(sharp(r))S(r)} implies
$$
0 = -\{r,r,q\} + (2r^2-\Delta(r,r))q - 3\Delta(r,q)r + \big(2T(r)\Delta(r,q) 
 - (r^\#,q) + \Delta(r_\# q,r) + N(r,q)\big)c.
$$

So, we have proved the identity~\eqref{triple(r,r,q)}.

The identity \eqref{triple(r,s,q)+triple(s,r,q)} is a linearization of \eqref{triple(r,r,q)},
while~\eqref{(r,s,q)} is a consequence of~\eqref{triple(r,s,q)+triple(s,r,q)}.
\hfill $\square$

\section{$\Psi$-map}

Define $\Psi(r,s,q)$ as follows,
\begin{multline} \label{PsiDef}
\Psi(r,s,q) = (r,s,q) - \Delta(q,s)r + \Delta(r,s)q \\
 + 1/4(\Delta(q_\#s,r) - \Delta(r_\#s,q) + 2T(r)\Delta(s,q) 
 - 2T(q)\Delta(r,s) + (q_\# s,r) - (r_\# s, q))c. \!\!\!  
\end{multline} 
By the definition, 
$\Psi(r,s,q) + \Psi(q,s,r) = 0$ for all $r,s,q$ and 
$\Psi(r,s,q) = 0$ if either of $r,s,q$ equals to~$c$.

The equalities~\eqref{(r,s,q)} and~\eqref{PsiDef} joint imply
$$ 
4\Psi(r,s,q) 
 = (r_\# s)_\# q - (s_\# q)_\# r
+ ((r,s) + 3\Delta(r,s))q  - ((q,s) + 3\Delta(q,s))r. 
$$

We introduce the following notations:
$$
\widetilde{(r,q)} = (r,q) + 3\Delta(r,q), \quad
(r,s,q)_\# = (r_\# s)_\# q - r_\#(s_\# q). 
$$
Then the last relation obtained has the form
\begin{equation} \label{PsiViaSharpAsso}
4\Psi(r,s,q) 
= (r,s,q)_\# + \widetilde{(r,s)}q - \widetilde{(s,q)}r.
\end{equation}
Similarly, we derive that
$$
4\Psi(r,s,q) 
= U_{q,s}(r) - U_{r,s}(q) + 3(\Delta(r,s)q - \Delta(q,s)r).
$$

We may rewrite~\eqref{T(sharpproduct(r,q))} as follows,
\begin{equation}
T(r_\# q) = T(r)T(q) - \widetilde{(r,q)}. \label{T(sharp)TildeInner}
\end{equation}

{\bf Lemma 5}.
Given a vector space~$V$, let $(N,\Delta,\#,c)$ be a~generalized sharped cubic form on~$V$.
Then the following identities are equivalent: 
\begin{gather}
\widetilde{(r_\# s,q)} = \widetilde{(r,s_\# q)}, \label{InvFormUnderTildeInner} \\
\{r,r^\#,q\} = (2N(r)-\Delta(r,r^\#))q - 3\Delta(r^\#,q)r, \label{triple(r,sharp(r),q)} \\
T(\Psi(r,s,q)) = 0. \label{T(Psi)=0} 
\end{gather}

{\sc Proof}.
We deduce~\eqref{triple(r,sharp(r),q)}:
\begin{multline*}
\{r,r^\#,q\} 
\mathop{=}^{\eqref{tripleproduct}} (r,r^\#)q + (r^\#,q)r - {r^\#}_\#(r_\# q) \\
\mathop{=}^{\eqref{Adjoint'}} (r,r^\#)q + (r^\#,q)r
 - (N(r)+\Delta(r^\#,r))q - (N(r,q)+\Delta(r^\#,q)+\Delta(r,r_\#q))r \\
 \mathop{=}^{\eqref{def:GSharp1},\,\eqref{(sharp(r),r)}} (2N(r)-\Delta(r,r^\#))q 
 - (-2/3(r^\#,q)+1/3(r_\#q,r) +\Delta(r^\#,q)+\Delta(r,r_\#q) )r,
\end{multline*}
which is equal to the right-hand side of~\eqref{triple(r,sharp(r),q)}
if and only~\eqref{InvFormUnderTildeInner} holds for $s = r$.
To show that~\eqref{InvFormUnderTildeInner} and~\eqref{triple(r,sharp(r),q)} 
are equivalent, it remains to derive~\eqref{InvFormUnderTildeInner} from itself fulfilled for $s = r$.
A linearization of
\begin{equation} \label{InvFormUnderTildeInnerPartial}
(r_\#q,r) + 3\Delta(r_\#q,r)
= (r_\#r,q) + 3\Delta(r_\#r,q)
\end{equation}
gives
$$
(r_\#q,s) + (s_\#q,r) + 3\Delta(r_\#q,s) + 3\Delta(s_\#q,r)
= 2(r_\#s,q) + 6\Delta(r_\#s,q).
$$
We may rewrite the last expression with the help of~\eqref{InvFormUnderSharpCycle}:
$$
N(r,s,q) + \Delta(r_\#q,s) + \Delta(s_\#q,r) + \Delta(r_\#s,q)
= (r_\#s,q) + 3\Delta(r_\#s,q).
$$
Because of the symmetry,~\eqref{InvFormUnderTildeInner} follows.

Due to~\eqref{PsiDef} and to Corollary 1, we have
\begin{multline*}
4T(\Psi(r,s,q)) 
 = 4(rs,q) - 4(r,sq) - 4T(r)\Delta(q,s) + 4T(q)\Delta(r,s) \\
 + 3(\Delta(q_\#s,r) - \Delta(r_\#s,q) + 2T(r)\Delta(q,s) - 2T(q)\Delta(r,s)
 + (r,s_\#q) - (r_\# s, q)) \\
 \mathop{=}\limits^{\eqref{product}} 2(r_\#s,q) + 2T(r)(s,q) 
 - 2S(r,s)T(q) - 2(r,s_\# q) - 2T(q)(r,s) + 2S(q,s)T(r) \\
 + 2T(r)\Delta(q,s) - 2T(q)\Delta(r,s) + 3(\Delta(q_\#s,r) 
 - \Delta(r_\#s,q) + (r,s_\#q) - (r_\# s, q)) \\
 \mathop{=}\limits^{\eqref{GNorm}} (r,s_\# q) - (r_\#s,q) + 2T(r)T(s)T(q) 
 - 2T(s)T(r)T(q) + 3(\Delta(q_\#s,r) - \Delta(r_\#s,q)) \\
 = \widetilde{(r,s_\# q)} - \widetilde{(r_\# s,q).}
\end{multline*}
Hence,~\eqref{InvFormUnderTildeInner} and~\eqref{T(Psi)=0} are equivalent. 
\hfill $\square$

{\bf Remark 3}.
Given a vector space~$V$, let $(N,\Delta,\#,c)$ be a~generalized sharped cubic form on~$V$ 
such that the form $(\cdot,\cdot)$ is invariant, i.\,e.~\eqref{InvForm} holds.
Then the form $\widetilde{(\cdot,\cdot)}$ is $\#$-invariant if and only if the equality
\begin{equation} 
\Delta(r_\#q,s) - \Delta(r,q_\#s) = \frac{T(s)\Delta(r,q) - T(r)\Delta(q,s)}{3} 
 \label{InvFormUnderDelta}
\end{equation}
if fulfilled for all $r,s,q\in V$.
When $N$ is inner, then invariancy of $(\cdot,\cdot)$ implies~\eqref{InvFormUnderDelta} 
and so implies $\#$-invariancy of
$\widetilde{(\cdot,\cdot)}$.

{\bf Lemma 6}.
Given a vector space~$V$, let $(N,\Delta,\#,c)$ be a~generalized sharped cubic form on~$V$ 
such that the form $(\cdot,\cdot)$ is invariant.
Then the following identities are fulfilled: 
\begin{gather}
(U_r(q),s) = (q,U_r(s)) + T(s)\Delta(r^\#,q)-T(q)\Delta(r^\#,s), \label{UopInvariancy} 
\end{gather}

\vspace{-0.9cm}

\begin{multline} \label{Ur(Usharp(r))}
U_r(U_{r^\#}(q)) 
= [-3(r^\#,q)\Delta(r,r^\#) + (N(r)+\Delta(r,r^\#))(\Delta(r^\#,q)+\Delta(r,r_\# q) \\
+2/3(T(r)\Delta(r,q)-T(q)\Delta(r,r)))]r + (N(r)+\Delta(r,r^\#))^2q,
\end{multline}

{\sc Proof}.
The formula~\eqref{UopInvariancy} holds by~\eqref{InvFormUnderSharpprod}.

We rewrite with the help of~\eqref{def:GSharp2} and~\eqref{U_r(sharp(r))}:
\begin{multline*}
U_r(U_{r^\#}(q)) 
= U_r( (r^\#,q)r^\# - {(r^\#)^\#}_\# q ) \\
= (r^\#,q)(N(r)-2\Delta(r^\#,r))r - (N(r) + \Delta(r^\#,r))U_r(r_\# q).
\end{multline*}
Further, we apply~\eqref{Adjoint'} and~\eqref{InvFormUnderSharpprod},
\begin{multline*}
U_r(r_\# q)
= (r,r_\# q)r - {r^\#}_\#(r_\# q) 
= ((q,r_\# r) + T(q)\Delta(r,r) - T(r)\Delta(r,q))r \\
- (N(r)+\Delta(r^\#,r))q - (N(r,q)+\Delta(r^\#,q)+\Delta(r,r_\#q))r.
\end{multline*}
Thus,
$$
U_r(U_{r^\#}(q)) 
= Ar + (N(r) + \Delta(r^\#,r))^2 q,
$$
where again by~\eqref{InvFormUnderSharpprod} we reduce
\begin{multline*}
A = (r^\#,q)(N(r)-2\Delta(r^\#,r)) 
 - (N(r) + \Delta(r^\#,r))( 2(r^\#,q) + T(q)\Delta(r,r) - T(r)\Delta(r,q) 
 \allowdisplaybreaks  \\
 - N(r,q) - \Delta(r^\#,q) - \Delta(r,r_\#q) ) 
 = -3(r^\#,q)\Delta(r^\#,r) \\
 + (N(r) + \Delta(r^\#,r))(N(r,q) - (r^\#,q) + \Delta(r^\#,q)+\Delta(r,r_\# q) 
 + T(r)\Delta(r,q) - T(q)\Delta(r,r) ).
\end{multline*}
It remains to use~\eqref{def:GSharp1-old} to prove~\eqref{Ur(Usharp(r))}.
\hfill $\square$

Let us prove some further properties of $\Psi$.

{\bf Lemma 7}.
Given a vector space~$V$, let $(N,\Delta,\#,c)$ be a~generalized sharped cubic form on~$V$ 
such that $\widetilde{(\cdot,\cdot)}$ is $\#$-invariant. 
Then the following identities for $\Psi$ are fulfilled:
\begin{gather}
\Psi(r,s,q) + \Psi(s,q,r) + \Psi(q,r,s) = 0,  \label{PsiJacobi} \\
\widetilde{(\Psi(r,s,q),x)} + \widetilde{(\Psi(q,s,x),r)} + \widetilde{(\Psi(x,s,r),q)} = 0. 
 \label{Delta-Psi-Tilde}
\end{gather}
If, additionally, $(\cdot,\cdot)$ is invariant, then
\begin{equation}
\Delta(\Psi(r,s,q),x) + \Delta(\Psi(q,s,x),r) + \Delta(\Psi(x,s,r),q) = 0. \label{Delta-Psi} 
\end{equation}

{\sc Proof}.
The equality~\eqref{PsiJacobi} follows by~\eqref{PsiViaSharpAsso}.

Based on~\eqref{PsiViaSharpAsso}, we rewrite and get
\begin{multline*}
4(\widetilde{(\Psi(r,s,q),x)} + \widetilde{(\Psi(q,s,x),r)} + \widetilde{(\Psi(x,s,r),q)}) \\ 
= \widetilde{((r,s,q)_\#,x)} + \widetilde{((q,s,x)_\#,r)} + \widetilde{((x,s,r)_\#,q)} \\
+ \widetilde{(r,s)}\widetilde{(q,x)} - \widetilde{(s,q)}\widetilde{(r,x)}
+ \widetilde{(q,s)}\widetilde{(x,r)} - \widetilde{(s,x)}\widetilde{(q,r)}
+ \widetilde{(x,s)}\widetilde{(r,q)} - \widetilde{(s,r)}\widetilde{(x,q)} \\
\mathop{=}\limits^{\eqref{InvFormUnderTildeInner}}
\widetilde{(r_\# s,x_\# q)} - \widetilde{(s_\# q,x_\# r)} 
+ \widetilde{(q_\# s,r_\# x)} - \widetilde{(s_\# x,r_\# q)} 
+ \widetilde{(x_\# s,q_\# r)} - \widetilde{(s_\# r,q_\# x)}
= 0,    
\end{multline*}
as required. 

\newpage

Let us prove~\eqref{Delta-Psi}, for this, we write down
\begin{multline*}
4\Delta(\Psi(r,s,q),x)
\mathop{=}\limits^{\eqref{PsiViaSharpAsso}} 
\Delta((r_\#s)_\#q - r_\#(s_\# q) + \widetilde{(r,s)}q - \widetilde{(q,s)}r,x) \\
\mathop{=}\limits^{\eqref{InvFormUnderDelta}} 
\Delta(r_\# s,q_\#x) - \Delta(s_\# q,r_\#x) 
 + \widetilde{(r,s)}\Delta(q,x) - \widetilde{(q,s)}\Delta(r,x) 
 \allowdisplaybreaks \\
+ \frac{T(x)\Delta(r_\# s,q)-T(r_\#s)\Delta(q,x)-T(x)\Delta(r,s_\# q)+T(s_\#q)\Delta(r,x)}{3} \\
\mathop{=}\limits^{\eqref{InvFormUnderDelta}} 
\Delta(r_\# s,q_\#x) - \Delta(s_\#q,r_\#x) 
 + \widetilde{(r,s)}\Delta(q,x) - \widetilde{(q,s)}\Delta(r,x) \\
 + \frac{T(x)(T(q)\Delta(r,s)-T(r)\Delta(s,q))}{9} 
 + \frac{T(s_\#q)\Delta(r,x)-T(r_\#s)\Delta(q,x)}{3}. 
\end{multline*}
Analogously, we have
\begin{multline*}
4\Delta(\Psi(q,s,x),r) 
 = \Delta(q_\#s,r_\#x) - \Delta(s_\#x,r_\#q) 
 + \widetilde{(s,q)}\Delta(r,x) - \widetilde{(s,x)}\Delta(r,q) \\
 + \frac{T(r)(T(x)\Delta(s,q)-T(q)\Delta(s,x))}{9} 
 + \frac{T(s_\#x)\Delta(r,q)-T(q_\#s)\Delta(r,x)}{3},
\end{multline*}

\vspace{-0.7cm}

\begin{multline*}
4\Delta(\Psi(x,s,r),q) 
 = \Delta(x_\#s,r_\#q) - \Delta(s_\#r,q_\#x) + \widetilde{(x,s)}\Delta(r,q) 
 - \widetilde{(r,s)}\Delta(q,x) \\
 + \frac{T(q)(T(r)\Delta(s,x)-T(x)\Delta(r,s))}{9} 
 + \frac{T(s_\#r)\Delta(q,x)-T(x_\#s)\Delta(r,q)}{3}.
\end{multline*}
The sum of the three expressions equals~0.
\hfill $\square$

{\bf Lemma 8}.
Given a vector space~$V$, let $(N,\Delta,\#,c)$ be a~generalized sharped cubic form on~$V$. 
Suppose that $(\cdot,\cdot)$ is invariant and nondegenerate,
$\widetilde{(\cdot,\cdot)}$ is $\#$-invariant, and $\dim V\geq2$.
Then $(N,\Delta,\#,c)$ is inner if and only if 
\begin{equation} \label{Delta-Psi-Sharp}
\Delta(s,\Psi(r,s,q)_\# x + \Psi(q,s,x)_\#r + \Psi(x,s,r)_\#q) = 0    
\end{equation}
holds for all $r,s,q,x\in V$.

{\sc Proof}. 
Let us rewrite~\eqref{Delta-Psi-Sharp} in more convenient form.
With the help of~\eqref{PsiViaSharpAsso} and~\eqref{InvFormUnderDelta}, we get 
\begin{multline} \label{Delta-Psi-Sharp-help}
4\Delta(s,\Psi(r,s,q)_\# x + \Psi(q,s,x)_\#r + \Psi(x,s,r)_\#q)  
= \Delta(s,(r,s,q)_\# x + \widetilde{(r,s)}q_\# x - \widetilde{(q,s)}r_\# x \\
+ (q,s,x)_\#r + \widetilde{(q,s)}x_\# r - \widetilde{(x,s)}q_\# r
+ (x,s,r)_\#q +\widetilde{(x,s)}r_\# q - \widetilde{(r,s)}q_\#x ) \\
\allowdisplaybreaks 
= \Delta((r,s,q)_\#, x_\#s) + \frac{T(s)}{3}\Delta((r,s,q)_\#,x)
- \frac{T((r,s,q)_\#)}{3}\Delta(x,s)  \\
+ \Delta((q,s,x)_\#, r_\#s) + \frac{T(s)}{3}\Delta((q,s,x)_\#,r) 
- \frac{T((q,s,x)_\#)}{3}\Delta(r,s) \\
+ \Delta((x,s,r)_\#, q_\#s) + \frac{T(s)}{3}\Delta((x,s,r)_\#,q) 
- \frac{T((x,s,r)_\#)}{3}\Delta(q,s).
\end{multline}
The sum of the three summands at $T(s)/3$ is zero due to~\eqref{Delta-Psi}.
Further, 
\begin{multline*}
\Delta((r_\#s)_\#q,x_\#s) - \Delta((q_\#(s_\#x),r_\#s)
\mathop{=}\limits^{\eqref{InvFormUnderDelta}}
\frac{1}{3}( T(x_\#s)\Delta(r_\#s,q) - T(r_\#s)\Delta(x_\#s,q) ) \\
\mathop{=}\limits^{\eqref{T(sharp)TildeInner}}
\frac{1}{3}( (T(x)T(s) - \widetilde{(x,s)})\Delta(r_\#s,q) 
- (T(r)T(s) - \widetilde{(r,s)})\Delta(x_\#s,q) ). 
\end{multline*}
Hence,
\begin{multline*}
\Delta((r,s,q)_\#, x_\#s) 
+ \Delta((q,s,x)_\#, r_\#s)
+ \Delta((x,s,r)_\#, q_\#s) \\
=  \frac{1}{3}( (T(x)T(s) - \widetilde{(x,s)})(\Delta(r_\#s,q)-\Delta(r,s_\#q))
- (T(r)T(s) - \widetilde{(r,s)})(\Delta(x_\#s,q)-\Delta(x,s_\#q)  \\
+ (T(q)T(s) - \widetilde{(q,s)})(\Delta(x_\#s,r)-\Delta(x,s_\#r)) ) \allowdisplaybreaks \\
\mathop{=}\limits^{\eqref{InvFormUnderDelta}}
\frac{1}{9}(
(T(x)T(s) - \widetilde{(x,s)})(T(q)\Delta(r,s)-T(r)\Delta(q,s)) \\
- (T(r)T(s) - \widetilde{(r,s)})(T(q)\Delta(x,s)-T(x)\Delta(s,q))  \\
+ (T(q)T(s) - \widetilde{(q,s)})(T(r)\Delta(x,s)-T(x)\Delta(r,s))
),
\end{multline*}
where the last expression equals the following one
\begin{multline} \label{Delta-Psi-Sharp-Equi}
T(q)\Delta(s,x)(r,s)-T(r)\Delta(s,x)(s,q)+T(x)\Delta(r,s)(s,q)-T(q)\Delta(r,s)(s,x)\\
+ T(r)\Delta(s,q)(s,x)-T(x)\Delta(s,q)(r,s) = 0
\end{multline}
with coefficient $1/9$.

The rest summands of~\eqref{Delta-Psi-Sharp-help} give by~\eqref{PsiViaSharpAsso} and~\eqref{T(Psi)=0}:
\begin{multline*}
\frac{1}{3}(
(T(q)\widetilde{(r,s)} - T(r)\widetilde{(q,s)})\Delta(x,s)
+ (T(x)\widetilde{(q,s)} - T(q)\widetilde{(x,s)})\Delta(r,s) \\
+ (T(r)\widetilde{(x,s)} - T(x)\widetilde{(r,s)})\Delta(q,s) ),
\end{multline*}
which is equal to~\eqref{Delta-Psi-Sharp-Equi} with coefficient $1/3$.

Therefore, we have showed that~\eqref{Delta-Psi-Sharp} is equivalent to~\eqref{Delta-Psi-Sharp-Equi}.
It is easy to check that if $(N,\Delta,\#,c)$ is inner, then~\eqref{Delta-Psi-Sharp-Equi} holds. 

Now, we want to show that if \eqref{Delta-Psi-Sharp-Equi} is true, then
the cubic form is inner.
Putting $q = c$ in~\eqref{Delta-Psi-Sharp-Equi}, we derive
\begin{equation} \label{CubicFormInnerEqui}
\Delta(s,x)\bigg((r,s)-\frac{T(r)T(s)}{3}\bigg) 
= \Delta(s,r)\bigg((x,s)-\frac{T(x)T(s)}{3}\bigg).
\end{equation}
Analogously, we write down
\begin{equation} \label{CubicFormInnerEqui'}
\Delta(s,r)\bigg((r,t)-\frac{T(r)T(t)}{3}\bigg) 
= \Delta(r,t)\bigg((r,s)-\frac{T(r)T(s)}{3}\bigg).
\end{equation}

Multiplying~\eqref{CubicFormInnerEqui} by~$\Delta(r,t)$ and adding~\eqref{CubicFormInnerEqui'} 
multiplied by~$\Delta(s,x)$, we get
$$
\Delta(r,s)\bigg(   
\Delta(s,x)\bigg((r,t)-\frac{T(r)T(t)}{3}\bigg) 
- \Delta(r,t)\bigg((x,s)-\frac{T(x)T(s)}{3}\bigg)
\bigg) = 0.
$$

If $\Delta\equiv0$, then $(N,\Delta,\#,c)$ is inner.

Otherwise, take $r,s$ such that $\Delta(r,s)\neq0$.
We may assume that $T(r) = 0$, since $\Delta(c,s) = 0$.
Now, we find $t\in V$ with the property $(r,t)\neq0$. 
Denote $\lambda = \Delta(r,t)/(r,t)$.
Hence, $\Delta(s,x) = \lambda \big((x,s)-\frac{T(x)T(s)}{3}\big)$
for all $x$ and all $s$~satisfying $\Delta(r,s)\neq0$ with fixed~$r$.
In particular, $\Delta(s,r) = \lambda \big((r,s)-\frac{T(r)T(s)}{3}\big)$.
Consider~$s$ such that $\Delta(r,s) = 0$.
Then by~\eqref{CubicFormInnerEqui'}, 
$\Delta(r,p)(r,s) = 0$ for all $p$. Hence, $(r,s) = 0$ and again 
\begin{equation} \label{criterion-inner}
\Delta(r,s) = \lambda ((r,s)-T(r)T(s)/3) 
\end{equation}
holds.

If for every $r\neq0$ such that $T(r) = 0$, one may find a corresponding~$s$ 
with the property $\Delta(r,s)\neq0$, then we have~\eqref{criterion-inner}.
Hence, every $a \in V$ may be written as $\mu c + r$, where $\mu\in F$ and $T(r) = 0$. 
Then
$\Delta(a,b) = \lambda \big((a,b)-\frac{T(a)T(b)}{3}\big)$,
where $\lambda\neq0$~depends on~$a$ and some~$t$.
From~\eqref{CubicFormInnerEqui}, we conclude that $\lambda$ is a~constant.

If we may find $r\neq0$ such that $T(r) = 0$ and $\Delta(r,s) = 0$ for all $s\in V$, 
then by~\eqref{CubicFormInnerEqui}, $\Delta(s,x) = 0$ for all $x$ and $s\not\perp r$ 
with respect to the form~$(\cdot,\cdot)$. Let us take any $a$ orthogonal to $r$ and fixed $s\not\perp r$. 
Then 
$\Delta(a,x) = \Delta(a+s,x) - \Delta(s,x) = 0$,
so $\Delta \equiv 0$, a contradiction.
Thus, $(N,\Delta,\#,c)$ is inner.
\hfill $\square$

Recall that the generalized sharped cubic form on the split spin factor $S(\alpha,E)$ is 
inner with $\lambda = \dfrac{3\alpha(1-\alpha)}{(1+\alpha)(\alpha-2)}$.
Therefore, the relations~\eqref{Delta-Psi} and~\eqref{Delta-Psi-Sharp} are fulfilled on $S(\alpha,E)$. 
In~\S6, we will prove that these identities hold on $S(\alpha,t,E)$.

{\bf Corollary 2}.
Given a vector space~$V$, 
let $(N,\Delta,\#,c)$ be an inner~generalized sharped cubic form on~$V$. 
Suppose that $(\cdot,\cdot)$ is invariant and nondegenerate,
$\widetilde{(\cdot,\cdot)}$ is $\#$-invariant, and $\dim V\geq2$.
Then 

a) $\Delta(s,x)(r,s) = \Delta(r,s)(s,x)$,
when either $T(s) = 0$ or $T(r) = T(x) = 0$,

b) $\Delta(s,\Psi(r,s,q)) = 0$ for all $r,s,q\in V$.

{\sc Proof}. 
a) The relation~\eqref{CubicFormInnerEqui} is equivalent by~\eqref{T,SOnUnit} to
$\Delta(s,x)(r,s) = \Delta(r,s)(s,x)$,
when either $T(s) = 0$ or $T(r) = T(x) = 0$.

b) We consider~\eqref{Delta-Psi-Sharp} with $x = c$.
Since $\Psi(q,s,c) = \Psi(c,s,r) = 0$ and by~\eqref{T(Psi)=0}, we get $\Delta(s,\Psi(r,s,q)) = 0$.
\hfill $\square$

\section{Sum of the three associators identity}

Now, we are ready to prove that $S(\alpha,t,E)$ satisfies the identity
\begin{equation} \label{Wb-identity}
W_b(a,c,d) := ((a,b,c),d,b) + ((c,b,d),a,b) + ((d,b,a),c,b) = 0.
\end{equation}

First, we show that $\widetilde{(\cdot,\cdot)}$ is $\#$-invariant on $S(\alpha,t,E)$.
By the proof of Lemma~5, it is enough to check~\eqref{InvFormUnderTildeInnerPartial}.
We compute for $r = az_1 + bz_2 + v$ and $q = gz_1 + hz_2 + w$:
\begin{multline*}
(r_\#r,q) + 3\Delta(r_\#r,q)    
= 2(1+\alpha)g( (\alpha a +\bar{\alpha}b )b + (\alpha-1)(t-1)\langle v,v\rangle) \\
+ 2(2-\alpha)h( (\alpha a +\bar{\alpha}b )a + \alpha(t-1)\langle v,v\rangle) 
- 2(1+\alpha+(2-\alpha)t)(\bar{\alpha}a+\alpha b)\langle v,w\rangle \\
+ 6\alpha(\alpha-1)(g-h)( (\alpha a + \bar{\alpha}b)(b-a)-(t-1)\langle v,v\rangle )
+ 6(\bar{\alpha}+\alpha t)(\bar{\alpha}a+\alpha b)\langle v,w\rangle;
\end{multline*}

\vspace{-0.7cm}

\begin{multline*}
(r_\#q,r) + 3\Delta(r_\#q,r)
= (1+\alpha)a(h(\alpha a+\bar{\alpha}b) + b(\alpha g+\bar{\alpha}h) 
 + 2(\alpha-1)(t-1)\langle v,w\rangle) \\
 + (2-\alpha)b(g(\alpha a+\bar{\alpha}b) + a(\alpha g+\bar{\alpha}h) 
 + 2\alpha(t-1)\langle v,w\rangle) \\
 - (1+\alpha +(2-\alpha)t)( (\bar{\alpha}a+\alpha b)\langle v,w\rangle 
 + (\bar{\alpha}g+\alpha h)\langle v,v\rangle ) 
 \allowdisplaybreaks \\
 + 3\alpha(\alpha-1)(a-b)( (\alpha a + \bar{\alpha}b)(h-g) + 
 (\alpha g + \bar{\alpha}h)(b-a) - 2(t-1)\langle v,w\rangle) \\
 + 3(\bar{\alpha}+\alpha t)( (\bar{\alpha}a+\alpha b)\langle v,w\rangle 
 + (\bar{\alpha}g+\alpha h)\langle v,v\rangle ).
\end{multline*}
In both $\widetilde{(r_\#r,q)}$ and $\widetilde{(r_\#q,r)}$, 
the coefficients at $\langle v,v\rangle$ equal
$2(2\alpha-1)(t-1)(\bar{\alpha}g+\alpha h)$, at $\langle v,w\rangle$ equal 
$4(2\alpha-1)(t-1)(\bar{\alpha}a+\alpha b)$. 
The rest summands equal to the same expression
$$
2(\alpha a+\bar{\alpha}b)( (1+\alpha)bg + (2-\alpha)ah + 3\alpha(\alpha-1)(a-b)(h-g)).
$$

Denote the coefficient at $c$ in~\eqref{PsiDef} as $\Phi(r,s,q)$.
If $\widetilde{(\cdot,\cdot)}$ is $\#$-invariant, then
\begin{equation}\label{PhiSimple}
\Phi(r,s,q) = 1/2(T(r)\Delta(s,q) - T(q)\Delta(r,s) 
+ \Delta(r_\#s,q) - \Delta(q_\#s,r) ). 
\end{equation}
Below, we apply that $\widetilde{(\cdot,\cdot)}$ is $\#$-invariant on $S(\alpha,t,E)$.

By~\eqref{PsiDef},
\begin{multline*}
((r,s,q),x,s) + ((q,s,x),r,s) + ((x,s,r),q,s) \\
 = \Delta(s,q)(r,x,s) - \Delta(r,s)(q,x,s)  + (\Psi(r,s,q),x,s)
 + \Delta(s,x)(q,r,s) - \Delta(q,s)(x,r,s) \allowdisplaybreaks \\ 
 + (\Psi(q,s,x),r,s) + \Delta(r,s)(x,q,s) - \Delta(s,x)(r,q,s)  + (\Psi(x,s,r),q,s) 
 \allowdisplaybreaks \\
 = \Delta(s,q)((r,s,x)+\Psi(x,s,r)) + \Delta(r,s)((x,s,q)+\Psi(q,s,x)) + \Delta(s,x)((q,s,r)+\Psi(r,s,q)) \\
 - ( \Delta(\Psi(r,s,q),x) + \Delta(\Psi(q,s,x),r) + \Delta(\Psi(x,s,r),q) )s \\
 - ( \Phi(\Psi(r,s,q),x,s) + \Phi(\Psi(q,s,x),r,s) + \Phi(\Psi(x,s,r),q,s) )c + \Psi_0.
\end{multline*}
where 
\begin{equation} \label{Psi0}
\Psi_0 = \Psi(\Psi(r,s,q),x,s) + \Psi(\Psi(q,s,x),r,s) + \Psi(\Psi(x,s,r),q,s).
\end{equation}

With the help of~\eqref{PsiDef} and Lemma~5, we rewrite the last expression as follows,
\begin{multline*}
((r,s,q),x,s) + ((q,s,x),r,s) + ((x,s,r),q,s) \\
 = \Delta(s,q)\left(-\Delta(r,s)x+\Delta(x,s)r + \frac{T(x)\Delta(r,s) 
 	- T(r)\Delta(x,s) + \Delta(x_\#s,r) - \Delta(r_\#s,x)}{2}c\right) \\
 + \Delta(r,s)\left(-\Delta(x,s)q+\Delta(q,s)x + \frac{T(q)\Delta(s,x) 
 	- T(x)\Delta(s,q) + \Delta(q_\#s,x) - \Delta(x_\#s,q)}{2}c\right) \\
 + \Delta(s,x)\left(-\Delta(s,q)r+\Delta(r,s)q + \frac{T(r)\Delta(s,q) 
 	- T(q)\Delta(r,s) + \Delta(r_\#s,q) - \Delta(q_\#s,r)}{2}c\right) \\
 - ( \Delta(\Psi(r,s,q),x) + \Delta(\Psi(q,s,x),r) + \Delta(\Psi(x,s,r),q) )s \\
 - ( \Phi(\Psi(r,s,q),x,s) + \Phi(\Psi(q,s,x),r,s) + \Phi(\Psi(x,s,r),q,s) )c + \Psi_0 \\
 = \frac{1}{2}(  \Delta(s,q)(\Delta(x_\#s,r) - \Delta(r_\#s,x))
 + \Delta(r,s)( \Delta(q_\#s,x) - \Delta(x_\#s,q) ) \\
 + \Delta(s,x)( \Delta(r_\#s,q) - \Delta(q_\#s,r) ) 
 - \frac{1}{2}\Delta(s,\Psi(r,s,q)_\# x + \Psi(q,s,x)_\#r + \Psi(x,s,r)_\#q) \\
 \allowdisplaybreaks
 + \frac{1}{2}( \Delta(\Psi(r,s,q),x_\#s) 
 + \Delta(\Psi(q,s,x),r_\#s) + \Delta(\Psi(x,s,r),q_\#s) ) \\
 + ( \Delta(\Psi(r,s,q),x) + \Delta(\Psi(q,s,x),r) + \Delta(\Psi(x,s,r),q) )((1/2)T(s)c-s)
 + \Psi_0.
\end{multline*}
Further, we will show that~$\Psi_0 = 0$, the identities~\eqref{Delta-Psi} 
and~\eqref{Delta-Psi-Sharp} are fulfilled on $S(\alpha,t,E)$ as well as 
\begin{equation} \label{Wb:additional}
\Delta(\Psi(r,s,q),x_\#s) + \Delta(\Psi(q,s,x),r_\#s) + \Delta(\Psi(x,s,r),q_\#s) = 0,    
\end{equation}

\vspace{-0.75cm}

\begin{multline} \label{Long-Delta-OnS}
\Delta(s,q)(\Delta(x_\#s,r) - \Delta(r_\#s,x))
+ \Delta(r,s)( \Delta(q_\#s,x) - \Delta(x_\#s,q) ) \\
+ \Delta(s,x)( \Delta(r_\#s,q) - \Delta(q_\#s,r) ) = 0.
\end{multline}

Now, let us explain that the identity~\eqref{Wb-identity} does not hold in general 
even for inner cubic forms. Consider the trivial case $\Delta \equiv 0$, 
which may be interpreted as an inner case with $\lambda = 0$.
Then the product coming from a sharped cubic form $(N,\#,c)$ is known to be Jordan~\cite{McCrimmon}. 
To check if~\eqref{Wb-identity} is fulfilled for the Jordan algebra, it is enough to study 
the case of a~special Jordan algebra, since the identity has the degree five.

Let $J$ be a special Jordan algebra, i.\,e. $J$ is a subalgebra of $A^{(+)}$, where 
$A$ is an associative algebra and the product $\circ$ in $A^{(+)}$ 
is defined as follows, $a\circ b = ab + ba$. 
Then $(a,b,c)_\circ = bac - bca + cab - acb$. Further,
\begin{multline*}
((a,b,c)_\circ,d,b)_\circ
= b(a,b,c)_\circ d + d(a,b,c)_\circ b - db(a,b,c)_\circ - (a,b,c)_\circ bd \\
\allowdisplaybreaks
= b^2(acd-cad) + b(ca - ac)bd + db(ac - ca)b + d(ca-ac)b^2 \\
- db^2(ac-ca) + db(ca-ac)b - b(ac-ca)bd - (ca-ac)b^2 d \\
= b^2(acd-cad) + (dca-dac)b^2 - db^2(ac-ca) - (ca-ac)b^2 d.
\end{multline*}
Thus,
$$
((a,b,c)_\circ,d,b)_\circ 
{+} ((c,b,d)_\circ,a,b)_\circ
{+} ((c,b,d)_\circ,a,b)_\circ 
{=} b^2(acd-cad + cda-dca + dac-adc ) {+} {\ldots}, 
$$
where the first two letters of all rest summands differ from $b^2$. 
Hence, this expression is nonzero in the case of any associative algebra~$A$, 
which does not satisfy any identity of degree less than~6.
For example, the matrix algebra~$M_3(F)$ is such an algebra~\cite{Levitzki}.
The space $M_3(F)$ is equipped with the identity matrix as a basepoint, 
the determinant as a~norm, and a sharp map sends a matrix to its adjoint.
Then the associated algebra is isomorphic to $M_3(F)^{(+)}$.
Slightly different sharped cubic form on $H_3(C)$, the Hermitian matrices 
over a~Cayley---Dickson algebra $C$, 
defines the simple Jordan algebra of Albert type~\cite{McCrimmon}.

To derive the identity~\eqref{Wb-identity}, we need the following result.

{\bf Lemma 9}.
In $S(\alpha,t,E)$, we have
\begin{equation} \label{PsiFormula}
\Psi(r,s,q) 
= (2\alpha-1)(t-1)(\langle u,w\rangle v - \langle u,v\rangle w),
\end{equation}
where $r = r_0 + v$, $s = s_0 + u$, $q = q_0 + w$ f
or $r_0,s_0,q_0\in Fz_1+Fz_2$ and $v,u,w\in E$.

{\sc Proof}.
Let us express the associator of the elements 
$r = az_1 + bz_2 + v$, $s = kz_1 + lz_2 + u$ and $q = gz_1 + hz_2 + w$.
We compute $(rs)q$ applying~\eqref{product-derived}:
\begin{multline*}
(rs)q
= (ak + \langle v,u\rangle)gz_1 + (bl +t\langle v,u\rangle)hz_2 
+ ((\alpha k + \bar{\alpha}l)\langle v,w\rangle \\
+ (\alpha a + \bar{\alpha}b)\langle u,w\rangle)(z_1+tz_2)
+ (\alpha g + \bar{\alpha}h)(\alpha k + \bar{\alpha}l)v \\
+ (\alpha g + \bar{\alpha}h)(\alpha a + \bar{\alpha}b)u
+ (\alpha(ak + \langle v,u\rangle) + \bar{\alpha}(bl +t\langle v,u\rangle))w.
\end{multline*}
Analogously, we have
\begin{multline*}
(qs)r
= (gk + \langle w,u\rangle)az_1 + (hl +t\langle w,u\rangle)bz_2 
+ ((\alpha k + \bar{\alpha}l)\langle v,w\rangle \\
+ (\alpha g + \bar{\alpha}h)\langle u,v\rangle)(z_1+tz_2)
+ (\alpha a + \bar{\alpha}b)(\alpha k+\bar{\alpha}l)w \\
+ (\alpha g + \bar{\alpha}h)(\alpha a+\bar{\alpha}b)u
+ (\alpha(kg + \langle w,u\rangle) + \bar{\alpha}(lh +t\langle w,u\rangle))v.
\end{multline*}
Therefore,
\begin{multline*}
(r,s,q)
= (g\langle v,u\rangle - a\langle w,u\rangle)z_1 
+ t(h\langle v,u\rangle - b\langle w,u\rangle)z_2 \\
+ ( (\alpha a + \bar{\alpha}b)\langle u,w\rangle
- (\alpha g + \bar{\alpha}h)\langle u,v\rangle))(z_1+tz_2) \\
- ( \alpha(\alpha-1)(a-b)(k-l) - (\alpha+\bar{\alpha}t)\langle v,u\rangle )w
+ ( \alpha(\alpha-1)(k-l)(g-h) - (\alpha+\bar{\alpha}t)\langle w,u\rangle )v.
\end{multline*}

It remains to substitute all known summands in~\eqref{PhiSimple}:
\begin{multline*}
\Psi(r,s,q)
 = (r,s,q) - \Delta(q,s)r + \Delta(r,s)q \\
 + 1/2(T(r)\Delta(s,q) - T(q)\Delta(r,s) 
 + \Delta(r_\#s,q) - \Delta(q_\#s,r) )c \\
 = (g\langle v,u\rangle - a\langle w,u\rangle)z_1 
 + t(h\langle v,u\rangle - b\langle w,u\rangle)z_2 
 + ( (\alpha a + \bar{\alpha}b)\langle u,w\rangle
 - (\alpha g + \bar{\alpha}h)\langle u,v\rangle)(z_1+tz_2) \\
 - ( \alpha(\alpha-1)(a-b)(k-l) - (\alpha+\bar{\alpha}t)\langle v,u\rangle )w
 + ( \alpha(\alpha-1)(k-l)(g-h) - (\alpha+\bar{\alpha}t)\langle w,u\rangle )v 
 \allowdisplaybreaks \\
 + (\alpha(\alpha-1)(k-l)(g-h)-(\bar{\alpha}+\alpha t)\langle u,w\rangle)(((1+\alpha)a + (2-\alpha)b)(z_1+z_2)/2-r) \\
 - (\alpha(\alpha-1)(a-b)(k-l)-(\bar{\alpha}+\alpha t)\langle v,u\rangle)
 (((1+\alpha)g + (2-\alpha)h)(z_1+z_2)/2-q) \\
 + \big[ \alpha(\alpha-1)(g-h)\big( (\alpha a+\bar{\alpha}b)(l-k) 
 + (\alpha k+\bar{\alpha}l)(b-a) - 2(t-1)\langle v,u\rangle \big) \\
 - \alpha(\alpha-1)(a-b)\big( (\alpha g+\bar{\alpha}h)(l-k) 
 + (\alpha k+\bar{\alpha}l)(h-g) - 2(t-1)\langle u,w\rangle \big) \\
 + (\bar{\alpha} + \alpha t)( (\bar{\alpha}a + \alpha b)\langle u,w\rangle
 - (\bar{\alpha}g + \alpha h)\langle v,u\rangle) \big](z_1+z_2)/2.
\end{multline*}
At $z_1/2$ we have the coefficient
\begin{multline*}
2(g\langle v,u\rangle - a\langle w,u\rangle
 + (\alpha a + \bar{\alpha}b)\langle u,w\rangle
 - (\alpha g + \bar{\alpha}h)\langle v,u\rangle) \\
 + (\alpha(\alpha-1)(k-l)(g-h)
 -(\bar{\alpha}+\alpha t)\langle u,w\rangle)((-1+\alpha)a + (2-\alpha)b) \\
 \allowdisplaybreaks
 - (\alpha(\alpha-1)(a-b)(k-l)-(\bar{\alpha}+\alpha t)\langle v,u\rangle)
((-1+\alpha)g + (2-\alpha)h) \\
 + \alpha(\alpha-1)(g-h)\big( (\alpha a+\bar{\alpha}b)(l-k) 
 + (\alpha k+\bar{\alpha}l)(b-a) - 2(t-1)\langle v,u\rangle \big) \\
 - \alpha(\alpha-1)(a-b)\big( (\alpha g+\bar{\alpha}h)(l-k) 
 + (\alpha k+\bar{\alpha}l)(h-g) - 2(t-1)\langle u,w\rangle \big) \\
 + (\bar{\alpha} + \alpha t)( (\bar{\alpha}a + \alpha b)\langle u,w\rangle
 - (\bar{\alpha}g + \alpha h)\langle v,u\rangle).
\end{multline*}
At $\langle v,u\rangle$, we have 
\begin{multline*}
g(2-2\alpha+(-1+\alpha)(\bar{\alpha}+\alpha t)
 -2\alpha(\alpha-1)(t-1)-\bar{\alpha}(\bar{\alpha}+\alpha t)) \\
 + h(-2\bar{\alpha}+(2-\alpha)(\bar{\alpha}+\alpha t)+2\alpha(\alpha-1)(t-1)
 -\alpha(\bar{\alpha}+\alpha t))
 = 0.
\end{multline*}
Analogously, we have zero coefficient at $\langle w,u\rangle$.
The rest summands equal $\alpha(\alpha-1)$ multiplied by
\begin{multline*}
(k-l)(g-h)((-1+\alpha)a+(2-\alpha)b)
- (a-b)(k-l)((-1+\alpha)g+(2-\alpha)h) \\
+ (g-h)((\alpha a+\bar{\alpha}b)(l-k) + (\alpha k+\bar{\alpha}l)(b-a))
- (a-b)((\alpha g+\bar{\alpha}h)(l-k) + (\alpha k+\bar{\alpha}l)(h-g)),
\end{multline*}
which is zero.

\newpage

Analogously, we have zero coordinate at $z_2$.
Finally, we have
$$
\Psi(r,s,q)
= (-\alpha-\bar{\alpha}t+\bar{\alpha}+\alpha t)(\langle u,w\rangle v - \langle u,v\rangle w) 
= (2\alpha-1)(t-1)(\langle u,w\rangle v - \langle u,v\rangle w),
$$
as required.
\hfill $\square$

{\bf Remark 4}.
It is easy to clarify, why Corollary 2b is true in $S(\alpha,t,E)$.
Indeed, by~\eqref{PsiFormula}, we have
$\Delta(s,\Psi(r,s,q)) = 0$
for $r = r_0 + v$, $s = s_0 + u$, $q = q_0 + w$,
where $r_0,s_0,q_0\in Fz_1+Fz_2$, $v,u,w\in E$, since
$$
\langle u,\langle u,w\rangle v - \langle u,v\rangle w \rangle
= \langle u,v\rangle \langle u,w\rangle - \langle u,v \rangle\langle u,w\rangle 
= 0.
$$
Let $\mu = (2\alpha-1)(t-1)$. In the case $\dim E = 2$, 
take a basis $e_1,e_2$ of $E$ such that $\langle e_1,e_2\rangle = 0$.
Then for $v=v_1e_1+v_2e_2$, $u=u_1e_1+u_2e_2$, and $w=w_1e_1+w_2e_2$, we have 
$$
\Psi(r,s,q) = \mu(v_1w_2 - v_2w_1)(u_2e_1 - u_1e_2). 
$$
Thus, $\Psi(r,s,q)$ is proportional to the vector $u^\perp = u_2 e_1 - u_1 e_2$, 
which is orthogonal to $u$ with respect to $\langle \cdot,\cdot\rangle$.

{\bf Corollary 3}.
In $S(\alpha,t,E)$, the relation $N(\Psi(r,s,q)) = 0$ holds for all $r,s,q$.
Hence, $\Psi(r,s,q)^3 = S(\Psi(r,s,q))\Psi(r,s,q)$.
An algebra $A$, in which every element satisfies the equality $x^3 = \varphi(x,x)x$ 
for some bilinear form $\varphi$, is called pseudo-composition algebra~\cite{Meyberg}.

{\bf Corollary 4}.
In $S(\alpha,t,E)$, we have
$\Delta(\Psi(r,s,q),x_\#s) = \Delta(\Psi(r,s,q)_\#s,x)$.
Put $s_0 = kz_1+lz_2$.
Applying the definition and Remark~4, we get
\begin{gather*}
\Delta(\Psi(r,s,q)_\#s,x)
= -(\bar{\alpha}k+\alpha l)\Delta(\Psi(r,s,q),x)
= (\bar{\alpha}k+\alpha l)(\bar{\alpha}+\alpha t)\langle \Psi(r,s,q),y\rangle, \\
\Delta(\Psi(r,s,q),s_\#x)
= -(\bar{\alpha}+\alpha t)\langle \Psi(r,s,q),s_\#x|_{E}\rangle
= (\bar{\alpha}+\alpha t)(\bar{\alpha}k+\alpha l)\langle \Psi(r,s,q),y\rangle,
\end{gather*}
hence, the required formula is proved.

{\bf Remark 5}.
Let us fix $u\in E$, then the product $[v,w] := \Psi (v,u,w)$ is a Lie one~\cite{Svinolupov}, 
thus,
$$
\Psi(\Psi(v,u,w),u,x) + \Psi(\Psi(w,u,x),u,v) + \Psi(\Psi(x,u,v),u,w) = 0
$$
holds for all $v,u,w,x\in E$.
Further, the ternary product $[v,u,w]:= \Psi(v,w,u)$ defines a~Lie triple system, 
i.\,e. the following identities for $[\cdot,\cdot,\cdot]$ hold:
\begin{gather*}
[x,y,z] + [y,x,z] = 0, \quad
[x,y,z] + [y,z,x] + [z,x,y] = 0, \\
[x,y,[u,v,w]] = [[x,y,u],v,w] + [u,[x,y,v],w] + [u,v,[x,y,w]],
\end{gather*}
more about triple systems see~\cite{Jacobson}. 
We believe that such construction of a~Lie triple system via an inner product is known, 
however, we are not able to find a suitable reference.

{\bf Theorem 3}.
The identity~\eqref{Wb-identity} holds on the algebra $S(\alpha,t,E)$.

{\sc Proof}.
Denote $\mu = (2\alpha-1)(t-1)$.
The identity~\eqref{Delta-Psi} holds on $S(\alpha,t,E)$, since by Lemma~9,
for $r = r_0 + v$, $s = s_0 + u$, $q = q_0 + w$, and $x = x_0 + y$,
where $r_0,s_0,q_0,x_0\in Fz_1+Fz_2$, $v,u,w,y\in E$, we have
\begin{multline*}
\Delta(\Psi(r,s,q),x) + \Delta(\Psi(q,s,x),r) + \Delta(\Psi(x,s,r),q) \\
 = \mu( \langle u,w\rangle \Delta(v,x) - \langle u,v\rangle \Delta(w,x)
 + \langle u,y\rangle \Delta(w,r) - \langle u,w\rangle \Delta(y,r) \\
 + \langle u,v\rangle \Delta(y,q) - \langle u,y\rangle \Delta(v,q) ) 
 = 0.
\end{multline*}

Let us verify that~\eqref{Delta-Psi-Sharp} is fulfilled on~$S(\alpha,t,E)$.
Denote $s_0 = kz_1 + lz_2$.
We apply Lemma~9:
\begin{multline*}
\Psi(r,s,q)_\# x
 = \mu(\langle u,w\rangle v - \langle u,v\rangle w)_\# x
 = 2\mu(t-1)(-\bar{\alpha}z_1+\alpha z_2)( \langle u,w\rangle \langle v,y\rangle 
 - \langle u,v\rangle \langle w,y\rangle ) \\
 - \mu\chi(x)(\langle u,w\rangle v - \langle u,v\rangle w),
\end{multline*}
where $\chi(p_1 z_1+p_2 z_2 + \omega) = \bar{\alpha}p_1 + \alpha p_2$.
Hence, by Remark~4, we get
\begin{equation} \label{delta(s,sharpproduct(Psi,x))}
\Delta(s,\Psi(r,s,q)_\# x)
 = -2\mu(t-1)\alpha(\alpha-1)(k-l)( \langle u,w\rangle \langle v,y\rangle 
 - \langle u,v\rangle \langle w,y\rangle ).
\end{equation}
Define $\pi = -2\mu(t-1)\alpha(\alpha-1)(k-l)$.
Then
\begin{multline*}
\Delta(s,\Psi(r,s,q)_\# x + \Psi(q,s,x)_\#r + \Psi(x,s,r)_\#q) \\
 = \pi( \langle u,w\rangle \langle v,y\rangle - \langle u,v\rangle \langle w,y\rangle 
 + \langle u,y\rangle \langle v,w\rangle - \langle u,w\rangle \langle v,y\rangle
 + \langle v,u\rangle \langle w,y\rangle - \langle u,y\rangle \langle v,w\rangle )
 = 0.    
\end{multline*}

Let us prove the relation~\eqref{Wb:additional}. By Lemma~9,
\begin{multline*}
\Delta(\Psi(r,s,q),x_\#s)
= \mu\Delta(\langle u,w\rangle v - \langle u,v\rangle w,-\chi(x)u-\chi(s)y) \\
= \mu\chi(s)(\bar{\alpha}+\alpha t)
(\langle u,w\rangle \langle v,y\rangle - \langle u,v\rangle \langle w,y\rangle).
\end{multline*}
As above (see~\eqref{delta(s,sharpproduct(Psi,x))}), we conclude that~\eqref{Wb:additional} holds.

Now, we check~\eqref{Long-Delta-OnS}.
First, we involve~\eqref{InvFormUnderTildeInner} and then~\eqref{product} and~\eqref{GNorm}:
\begin{multline} \label{DeltaDelta}
- 3( \Delta(s,q)(\Delta(x_\#s,r) - \Delta(r_\#s,x))
 + \Delta(r,s)( \Delta(q_\#s,x) - \Delta(x_\#s,q) ) \\
 + \Delta(s,x)( \Delta(r_\#s,q) - \Delta(q_\#s,r) ) 
 = \Delta(s,q)( (x_\#s,r) - (r_\#s,x)) \\
 + \Delta(r,s)( (q_\#s,x) - (x_\#s,q) ) 
 + \Delta(s,x)( (r_\#s,q) - (q_\#s,r) \\
 = 2( \Delta(s,q)((xs,r)-(x,sr))
 + \Delta(r,s)( (qs,x) - (q,sx))
 + \Delta(s,x)( (rs,q) - (r,sq)) ) \\
 + \Delta(s,q)( T(x)\Delta(s,r) - T(r)\Delta(s,x) )
 + \Delta(s,r)( T(q)\Delta(s,x) - T(x)\Delta(s,q) ) \\
 + \Delta(s,x)( T(r)\Delta(s,q) - T(q)\Delta(s,r) ) \\
 = 2( \Delta(s,q)((xs,r){-}(x,sr))
 + \Delta(r,s)( (qs,x) {-} (q,sx))
 + \Delta(s,x)( (rs,q) {-} (r,sq)) ).
\end{multline}

By~\eqref{inner-product-derived} and~\eqref{product-derived}, we compute
\begin{multline} \label{invariancy-deduced}
(rs,q) - (r,sq)
= (1+\alpha)g(ak+\langle v,u\rangle) + (2-\alpha)h(bl+t\langle v,u\rangle)
\\
+ (1+\alpha+(2-\alpha)t)( (\alpha a+\bar{\alpha}b)\langle u,w\rangle 
+ (\alpha k+\bar{\alpha}l)\langle v,w\rangle  ) 
\allowdisplaybreaks \\
- (1+\alpha)a(gk+\langle u,w\rangle) + (2-\alpha)b(hl+t\langle u,w\rangle)
\\ 
+ (1+\alpha+(2-\alpha)t)( (\alpha g+\bar{\alpha}h\langle u,v\rangle 
+ (\alpha k+\bar{\alpha}l)\langle v,w\rangle  ) \\
= (1-\alpha^2+\alpha(\alpha-2)t)( (g-h)\langle v,u\rangle - (a-b)\langle u,w\rangle ).
\end{multline}
Denote $\nu = (1-\alpha^2+\alpha(\alpha-2)t)$ and let $x_0 = mz_1 + nz_2$. Thus, 
\begin{multline*}
\Delta(s,x)( (rs,q) - (r,sq))
= \nu\alpha(\alpha-1)(k-l)\big(
(m-n)(g-h)\langle v,u\rangle - (m-n)(a-b)\langle u,w\rangle
\big) \\
+ \nu(\bar{\alpha}+\alpha t)( 
(a-b)\langle u,w\rangle \langle u,y\rangle 
- (g-h)\langle u,v\rangle \langle u,y\rangle ).
\end{multline*}
The analogous expressions for 
$\Delta(s,q)((xs,r)-(x,sr))$ 
and $\Delta(r,s)( (qs,x) - (q,sx))$ joint provide that~\eqref{DeltaDelta} equals zero.

Finally, it remains to prove that $\Psi_0 = 0$, see~\eqref{Psi0}.
By Lemma~9, we express
$$
\Psi(\Psi(r,s,q),x,s)
 = \mu^2( \langle v,u\rangle \langle w,y\rangle u - \langle w,u\rangle \langle y,v\rangle u 
 - \langle y,u\rangle \langle v,u\rangle w + \langle w,u\rangle \langle y,u\rangle v),
$$
Hence, 
\begin{multline*}
\Psi_0
 = \mu^2(\langle v,u\rangle \langle w,y\rangle u - \langle w,u\rangle \langle y,v\rangle s 
 - \langle y,u\rangle \langle v,u\rangle w + \langle w,u\rangle \langle y,u\rangle v \\
 + \langle w,u\rangle \langle y,v\rangle u - \langle y,u\rangle \langle v,w\rangle u 
 - \langle v,u\rangle \langle w,u\rangle y + \langle y,u\rangle \langle v,u\rangle w \\
 + \langle y,u\rangle \langle v,w\rangle u - \langle v,u\rangle \langle w,y\rangle u 
 - \langle w,u\rangle \langle y,u\rangle v + \langle v,u\rangle \langle w,u\rangle y )
 = 0.
\end{multline*}
The statement is proved. \hfill $\square$

{\bf Remark 6}.
Due to~\eqref{invariancy-deduced}, we see that the form $(\cdot,\cdot)$ 
is invariant on $S(\alpha,E)$, as it was noted in~\cite{McInroy}.

{\bf Remark 7}.
Let us return to Example~1.
We may introduce the bilinear form 
$\widetilde{(r,q)} = (r,q) + \Delta(r,q)$.
Then 
$$
T(r_\# q) = (1-4\lambda)T(r)T(q) - \widetilde{(r,q)}, \quad
\widetilde{(r_\# s,q)} = \widetilde{(r,s_\# q)}.
$$
Define the trilinear form~$\Psi$ by the formula~\eqref{PsiViaSharpAsso}. 
Denote $r = (a,b,c)$, $s = (i,j,k)$, and $q = (e,f,g)$.
Thus, we have
\begin{multline*}
\Psi(r,s,q)
= \lambda( j(-ag+ce+bg-cf) + k (-af+be - bg+ cf), \\
i(ag-ce-bg+cf) + k(af-be-ag+ce), \
i(af-be+bg-cf) + j(-af+be+ag-ce)).
\end{multline*}
Then the relations~\eqref{PsiJacobi},~\eqref{Delta-Psi},~\eqref{Delta-Psi-Sharp} 
and~\eqref{Psi0} are fulfilled.
Further, the identity of the three associators holds, 
and the ternary product $[v,u,w]:= \Psi(v,w,u)$ defines a~Lie triple system, see the code in GAP~\cite{Code}.
Moreover, the identity~\eqref{Wb-identity} is fulfilled on the space 
$V = A\otimes_F F^3\cong A^{\otimes3}$, where $A = F[\lambda]$.

\section{Identities}

In this section, we prove that the  algebra~$S(\alpha,t,E)$ does not satisfy 
any polynomial identity of degrees 3 and 4, and all identities of degree~5 
satisfied by $S(\alpha,E)$ follow from commutativity and the identity~\eqref{Wb-identity}.
In~1989, S.Yu. Vasilovsky found a basis of the $T$-ideal of identities fulfilled 
on the simple Jordan algebra of a nondegenerate form considered 
over a~field of characteristic~0~\cite{Vasilovsky}. One of them has the close form 
$(d,(a,b,c),b) + (a,(c,b,d),b) + (c,(d,b,a),b) = 0$.

In~\cite{Osborn}, it was proved that if a commutative (non-associative) unital algebra~$A$ 
over a field of characteristics not 2 or 3 satisfies an identity of degree 4 not implied 
by the commutative law, then $A$ satisfies at least one of the following three identities:
\begin{gather}
(x^2x)x = x^2x^2, \label{id-1} \\
2((yx)x)x + yx^3 = 3(yx^2)x, \label{id-2} \\
2(y^2x)x + 2(x^2y)y +(yx)(yx) = 2((yx)y)x+2((yx)x)y+y^2x^2. \label{id-3}
\end{gather} 

Then we have the following:

{\bf Lemma 10}.
Let $E$ has a dimension $n\ge 1$ and $\alpha,t\notin \{0,1\}$. 
Then every identity of degree no more than 4 in the algebra~$S(\alpha,t,E)$ 
follows from commutativity.

{\sc Proof}.
To prove the statement, it is enough to show that the algebra~$S(\alpha,t,E)$ 
does not satisfy the identities \eqref{id-1}--\eqref{id-3}.

First, let us show the identity \eqref{id-1} does not hold. 
Consider the left hand-side of \eqref{id-1} and set $x = e\in E$ 
such that $\langle e,e\rangle = 1$, then we have 
$$
(e^2e)e=((z_1+t z_2)e)e)
= (\alpha e + t (1-\alpha) e)e 
= (\alpha+t(1-\alpha))(z_1+t z_2).
$$
The right-hand side of~\eqref{id-1} for $x=e$ gives
$$
e^2e^2=(z_1+t z_2)(z_1+t z_2)
= z_1+t^2z_2.
$$

Since $t\neq 0,1$, we conclude that the identity \eqref{id-1} does not hold.

To show that \eqref{id-2} does not hold, it is enough to consider 
$x = e\in E$ such that $\langle e,e\rangle = 1$ and $y=z_1$. 
Then the left-hand side of~\eqref{id-2} equals 
\begin{equation}\label{id2subl ez_1}
2((z_1 e)e)e+z_1 e^3 
= 3\alpha(\alpha + t(1-\alpha))e,
\end{equation}
while the right-hand side of~\eqref{id-2} gives  
\begin{equation}\label{id2subr ez_1}
3(z_1e^2)e=3\alpha e.
\end{equation}
We see that the right-hand sides of \eqref{id2subl ez_1} and \eqref{id2subr ez_1} 
are equal if and only if $t=1$. By the conditions, $t\neq 1$ and therefore 
the identity \eqref{id-2} does not hold.

Now we consider \eqref{id-3}. Define  
$$
\phi(x,y)=2(y^2x)x + 2(x^2y)y +(yx)(yx) - 2((yx)y)x-2((yx)x)y-y^2x^2.
$$

Then $\phi(e,z_1)=(1-\alpha^2)z_1+t\alpha (2-\alpha)z_2 \neq 0$ 
for any $\alpha,t\not\in\{0,1\}$. 
Consequently, \eqref{id-3} is not an identity in the algebra~$S(\alpha,t,E)$.   
\hfill $\square$

In~\cite{Osborn}, the list of all irreducible relative to commutativity identities 
of degree five is given. There are exactly five such identities, 
and the fourth of them~\cite[eq.\,(15)]{Osborn} with $\delta_2 = -\delta_1\neq0$ 
is nothing more than~\eqref{Wb-identity} with one of the three variables $a,c,d$ 
equal to $b$, e.\,g., $d = b$.

The proof of the following theorem  is established through computations conducted 
with the assistance of software programs such as Wolfram Mathematica and Albert~\cite{Albert}. 

{\bf Theorem 4}.
Let $E$ has a dimension $n\ge 2$ and $\alpha\notin \{-1,0,1/2,1,2\}$. 
Every identity of degree no more than~5 in the algebra~$S(\alpha,E)$ 
over a~field of characteristics 0 is a~consequence of commutativity and~\eqref{Wb-identity}.

{\sc Proof}.
By Lemma~10, it remains to show that there are no identities in degree 5, 
which do not follow from commutativity and the identity \eqref{Wb-identity}.

Let $\mathcal{W}(X)$ denote a free algebra defined by identities 
of commutativity and~\eqref{Wb-identity}, which is generated by a set~$X$. 
Since we deal with a~field of characteristics~0, then every polynomial 
identity is equivalent to a set of multilinear identities \cite{ZSSS}.

Let $\mathcal{P}$ be a monomial basis of the multilinear part of degree~5 
of the free commutative algebra $\Com(X)$.
Then $\mathcal{P}$ consists of the 60~monomials of the type $(((**)*)*)*$, 
30~monomials of the type $((**)*)(**)$, and 15~monomials of the type $((**)(**))*$. 
Define the set
$$
\begin{array}{ccccc}
\mathcal{Z}
 = &\{((x_3 x_5) x_4) (x_1 x_2), & ((x_4 x_5) x_3) (x_1 x_2), 
&((x_2 x_5) x_4) (x_1 x_3), & ((x_4 x_5) x_2) (x_1 x_3), \\
&((x_2 x_5) x_3) (x_1 x_4), & ((x_3 x_5) x_2) (x_1 x_4), 
&(((x_1 x_5) x_4) x_3) x_2, & (((x_2 x_5) x_4) x_3) x_1, \\
&(((x_3 x_5) x_4) x_2) x_1, & (((x_4 x_5) x_3) x_2) x_1\}.
\end{array}
$$

To construct a monomial basis~$\mathcal{B}$ of the multilinear part 
of degree~5 of $\mathcal{W}(X)$, where 
$X=\{x_1,x_2,x_3,x_4,x_5\}$, we employ the software program Albert and obtain 95~basic monomials. 
We can represent the set of multilinear basic monomials as 
$\mathcal{B}=\mathcal{P} \setminus \mathcal{Z}$.

If there exists a multilinear polynomial identity of degree~5 fulfilled on $S(\alpha,t, E)$, 
which does not follow from commutativity and~\eqref{Wb-identity}, 
then it can be represented as a~linear combination of monomials from $\mathcal{B}$. 
Let us define a linear combination of elements in~$\mathcal{B}$ as 
$$
\psi(x_1,x_2,x_3,x_4,x_5)=\sum_{b_i\in\mathcal{B}}\lambda_i b_i.
$$ 
To establish the theorem, it is necessary to demonstrate the linear independence 
of monomials from $\mathcal{B}$. To achieve this, we use the Wolfram Mathematica software tool. 
A special code has been developed for calculating all substitutions, 
extracting homogeneous equations from them and solving them \cite{Code}.

Since $\dim E\geq2$, it is enough to show that all identities in degree 5 
fulfilled on $S(\alpha,t,E_0)$, where $\dim E_0 = 2$, follow from commutativity and~\eqref{Wb-identity}. 
Let us choose a~basis $e,f$ of $E_0$ such that 
$\langle e,e\rangle = \langle f,f\rangle = 1$ 
and $\langle e,f\rangle = 0$.
We use the function \texttt{Tuples}$[\{z_1,z_2,e,f\},5]$ to generate all 1024 possible permutations 
of length~5 using the basic elements $\{z_1,z_2,e,f\}$ and substitute them into 
$\psi(x_1,x_2,x_3,x_4,x_5)$. Employing the function \texttt{Union[]}, 
we express the obtained polynomials in terms of the coefficients $\lambda_i$ 
and the basic elements $\{z_1,z_2,e,f\}$. This yields a set of 635 polynomials.

Further, we express these polynomials by collecting coefficients at the basic elements 
$\{z_1,z_2,e,f\}$ and extract the coefficients corresponding to these elements with 
the functions \texttt{Collect[]} and \texttt{Coefficient[]}. 
By employing the function \texttt{Union[]} once more, we reduce the number of polynomials to 498. 
Then we consider the system of equations formed by setting all these polynomials equal to zero. 
This system of equations is expressed in the coefficients $\lambda_i$, where $i\in \{1,\ldots, 95\}$.

The only trivial solution that emerges is $\lambda_i=0$ for $\alpha\notin\{-1,0,\frac{1}{2},1,2\}$, 
where $i\in \{1,\ldots, 95\}$. This result demonstrates that the monomials involved 
in the linear combination $\psi(x_1,x_2,x_3,x_4,x_5)$ are linearly independent. 
This completes the proof.
\hfill $\square$

\textbf{Remark 8}. 
The above theorem is valid when $t = (\alpha^2-1)/\alpha(\alpha-2)$. 
However, it is essential to note that for a general value of $t$ 
that does not depend on $\alpha$, the algebra~$S(\alpha,t,E)$ can have 
an identity of degree 5 which does not follow from commutativity and~\eqref{Wb-identity}.

For example, for $t=5$ and $\alpha=11/4$, there is an identity 
\begin{multline*}
((c,a,e),b,d) + ((e,a,d),b,c) + ((d,a,c),b,e) \\
+ (c,b,a)[R_d,R_e] + (d,b,a)[R_e,R_c] + (e,b,a)[R_c,R_d] = 0,
\end{multline*}
which does not follow from commutativity and~\eqref{Wb-identity}. 
The validity of the identity can be checked using a program given in~\cite{Kadyrov} 
or requiring a program from the authors.

\section{Open problems}

We finish the work with several open problems concerned the subject.

\begin{itemize}
\item Find an identity, which does not follow from commutativity and is fulfilled 
on every algebra associated to a generalized cubic form.
\item Does the identity~\eqref{InvFormUnderDelta} follow from the definition 
of generalized sharped cubic form, or there exists a counterexample to it? 
\item Given a generalized sharped cubic form $(N,\Delta,\#,c)$, 
which satisfies~\eqref{InvFormUnderDelta}, is it true that $N(\Psi(r,s,q)) = 0$ for all $r,s,q$?
\item Find the basis of the $T$-ideal of identities fulfilled on $S(\alpha,t,E)$.
\end{itemize}

\section{Acknowledgments}

V. Gubarev is supported by Mathematical Center in Akademgorodok under 
agreement No. 075-15-2022-281 with the Ministry of Science 
and Higher Education of the Russian Federation.
A.S. Panasenko is supported by the Program of fundamental scientific researches 
of Russian Academy of Sciences, project FWNF-2022-0002.

The results of \S2 are supported by the Program of fundamental scientific researches 
of Russian Academy of Sciences, project FWNF-2022-0002.
The results of \S4--6 are supported by Mathematical Center in Akademgorodok 
under agreement No. 075-15-2022-281 with the Ministry of Science 
and Higher Education of the Russian Federation.

\noindent Vsevolod Gubarev \\
Novosibirsk State University \\
Pirogova str. 2, 630090 Novosibirsk, Russia \\
Sobolev Institute of Mathematics \\
Acad. Koptyug ave. 4, 630090 Novosibirsk, Russia \\
e-mail: wsewolod89@gmail.com

\smallskip

\noindent Farukh Mashurov \\
Suleyman Demirel University \\
Abylai Khan Street 1/1 \\ Kaskelen,  Kazakhstan \\
e-mail: f.mashurov@gmail.com

\smallskip

\noindent Alexander Panasenko \\
Sobolev Institute of Mathematics \\
Acad. Koptyug ave. 4, 630090 Novosibirsk, Russia \\
Novosibirsk State University \\
Pirogova str. 2, 630090 Novosibirsk, Russia \\
e-mail: a.panasenko@g.nsu.ru

\begin{thebibliography}{99}
\bibitem{Albert} Albert Version 4.0M6; 
https://web.osu.cz/\,${}_{\widetilde{}}$\,Zusmanovich/soft/albert.

\bibitem{Elduque} 
A. Elduque and S. Okubo, 
On algebras satisfying $x^2x^2 = N(x)x$, Math. Z. (2) {\bf 235} (2000), 275--314.

\bibitem{Code} 
https://github.com/GubarevV/Generalized-sharped-cubic-form.

\bibitem{HRS} 
J.I. Hall, F. Rehren and S. Shpectorov, 
Primitive axial algebras of Jordan type, J.~Algebra {\bf 437} (2015), 79--115.

\bibitem{Jacobson} 
N. Jacobson, 
Lie and Jordan triple systems, Am. J. Math. (1) {\bf 71} (1949), 149--170.

\bibitem{Kadyrov} 
S.~Kadyrov, F.~Mashurov, 
Unified computational approach to nilpotent algebra classification problems, 
Comm.  Math. (2) {\bf 29} (2021), 215--226.

\bibitem{Levitzki} 
J. Levitzki, 
A theorem on polynomial identities, Proc. of AMS, {\bf 1} (1950), 334--341.

\bibitem{McCrimmon2} 
K. McCrimmon, 
The Freudenthal-Springer-Tits constructions of exceptional Jordan
algebras, Trans. Am. Math. Soc. {\bf 139} (1969), 495--510.

\bibitem{McCrimmon} 
K. McCrimmon, 
A Taste of Jordan Algebras. N.-Y.: Springer-Verlag, 2004.

\bibitem{McInroy} 
J. McInroy and S. Shpectorov, 
Split spin factor algebras, J. Algebra {\bf 595} (2022), 380--397.

\bibitem{Meyberg} 
K. Meyberg and J.M. Osborn,
Pseudo-composition algebras, Math. Z. {\bf 214} (1993), 67--77.

\bibitem{Osborn} 
J.M. Osborn, 
Identities of Non-Associative Algebras, Can. J. Math. {\bf 17} (1965), 78--92.

\bibitem{Svinolupov}
S.I. Svinolupov, V.V. Sokolov, 
Vector-matrix generalizations of classical integrable equations, 
Theoret. and Math. Phys. (2) {\bf 100} (1994), 959--962.

\bibitem{Vasilovsky} 
S. Vasilovsky,
A finite basis for polynomial identities of the Jordan algebra of a~bilinear form,
Tr. Inst. Mat. (Novosib.) {\bf 16} (1989),
Issled. Teor. Kolets Algebre, 5--37. 
English translation: Siberian Adv. Math. (4) {\bf 1} (1991), 142--185.

\bibitem{Walcher} 
S. Walcher, 
On algebras of rank three, Commun. Algebra {\bf 27} (1999), 3401--3438.

\bibitem{Yabe} 
T. Yabe,
On the classification of 2-generated axial algebras of Majorana type,
J.~Algebra {\bf 619} (2023), 347--382.

\bibitem{ZSSS}
K.A. Zhevlakov, A.M. Slin'ko, I.P. Shestakov, and A.I. Shirshov,
Rings that are nearly associative. N.-Y.: Academic Press, 1982.
\end{thebibliography}
\end{document}